\documentclass[preprint,sort&compress]{elsarticle}
\usepackage{amssymb,amsthm,amsmath,commath,bbm,psfrag,setspace,color,array,booktabs,multirow,url,float}
\usepackage[FIGTOPCAP,nooneline]{subfigure}

\pdfoutput=1

\usepackage{amsmath,amsthm,amssymb}
\usepackage[dvipsnames]{xcolor}
\usepackage{xcolor}
\usepackage{color}

\usepackage[utf8x]{inputenc}

\usepackage{mathscinet}
\usepackage{verbatim}
\usepackage{enumerate}
\usepackage{mathrsfs}
\usepackage[colorlinks,citecolor=red]{hyperref}

\newtheorem{theorem}{Theorem}[section]
\newtheorem{lemma}{Lemma}[section]
\newtheorem{remark}{Remark}[section]

\newtheorem{proposition}{Proposition}[section]

\newcommand{\RR}{\mathbb{R}}

\newcommand{\seq}[1]{\left<#1\right>}

\providecommand{\abs}[1]{\lvert#1\rvert}
\providecommand{\norm}[1]{\lVert#1\rVert}

\begin{document}

\begin{frontmatter}

\title{A non-linear stable non-Gaussian process in fractional time}

\author{S. Solís (1) \corref{mycorrespondingauthor}}

\author{V. Vergara (2)}

\cortext[mycorrespondingauthor]{Corresponding author.\\(1) Departamento de Matem\'aticas, Facultad de Ciencias Naturales y Matemáticas, Escuela Superior Politécnica del Litoral,\\  Guayaquil, Ecuador; e-mail: ssolis@espol.edu.ec.\\(2) Departamento de Matem\'aticas, Facultad de Ciencias F\'isicas y Matem\'aticas, Universidad de Concepci\'on,\\  Concepci\'on, Chile; e-mail: vvergaraa@udec.cl.}

\begin{abstract}
A subdiffusion problem in which the diffusion term is related to a stable stochastic process is introduced. Linear models of these systems have been studied in a general way, but non-linear models require a more specific analysis. The model presented in this work corresponds to a non-linear evolution equation with fractional time derivative and a pseudo-differential operator acting on the spatial variable. This type of equations has a couple of fundamental solutions, whose estimates for the $L_p-$norm are found to obtain three main results concerning mild and global solutions. The existence and uniqueness of a mild solution is based on the conditions required in some parameters, one of which is the order of stability of the stochastic process. The existence and uniqueness of a global solution is found for the case of small initial conditions and another for non-negative initial conditions. The relationship between Fourier analysis and Markov processes, together with the theory of fixed points in Banach spaces, is particularly exploited. In addition, the present work includes the asymptotic behavior of global solutions as a linear combination of the fundamental solutions with $L_p-$decay. 
\end{abstract}


\end{frontmatter}

\section{Introduction}
\label{intro}
Let $\lambda\in \RR$, $\gamma>1$. We consider the following Cauchy problem 
\begin{equation}
\label{general}
\begin{split}
\partial^{\alpha}_t(u-u_0)(t,x)+\Psi_{\beta}(-i\nabla)u(t,x) & = \lambda \abs{u(t,x)}^{\gamma-1} u(t,x),\quad t>0, \; x\in\mathbb{R}^d,\\ u(t,x)|_{t=0} & = u_0(x),\quad x\in\mathbb{R}^d, 
\end{split}
\end{equation}
where $\partial^{\alpha}_t$ is the \textit{Riemann–Liouville fractional derivative} of order $\alpha\in(0,1)$ given by
\[
\partial^{\alpha}_t v=\dfrac{d}{dt}\displaystyle\int_0^t g_{1-\alpha}(t-s)v(s)ds=:\dfrac{d}{dt}(g_{1-\alpha}*v)(t),
\]
with $g_{\rho}(t):=\frac{1}{\Gamma(\rho)}t^{\rho-1}$ and the Euler Gamma function $\Gamma(\cdot)$. The function $u_0$ stands for the initial data in a certain Lebesgue space. The term $\Psi_{\beta}(-i\nabla)$ is a \textit{singular integral operator} of order $\beta\in(0,2)$ with symbol $\psi$, that is
\[
\Psi_{\beta}(-i\nabla) v(x)= \mathcal{F}^{-1}_{\xi\to x}[\psi(\xi) (\mathcal{F}v)(\xi)], \, v\in C_0^{\infty}(\RR^d),
\]
where $\mathcal{F}$ is the Fourier transform in $\RR^d$ and $\mathcal{F}^{-1}$ is its inverse. As usual, $C_0^{\infty}(\RR^d)$ denotes the space of test functions on $\RR^d$. The symbol $\psi$ is a measurable function on $\RR^d$ given by
\[
\psi(\xi)=\norm{\xi}^\beta \omega_\mu\left(\frac{\xi}{\norm{\xi}}\right),\quad \xi\in\mathbb{R}^d,
\]
where
\begin{equation}
\label{omega_mu}
\omega_\mu(\theta):=\displaystyle\int_{S^{d-1}}\abs{\theta\cdot\eta}^\beta\mu(d\eta),\quad\theta\in S^{d-1}.
\end{equation}
Here, $\mu(d\eta)$ is a centrally symmetric finite (non-negative) Borel measure defined on the unitary sphere $S^{d-1}$, called \textit{spectral measure}, and $\omega_\mu$ is a continuous function on $S^{d-1}$, see e.g. \cite[Section 1.8]{Kol19}. Whenever $\mu(d\eta)=\varrho(\eta)d\eta$, where $\varrho$ is a continuous function on $S^{d-1}$, we will refer to $\varrho$ as the \textit{density of $\mu$}. Some restrictions on the function $\varrho$ may be required for the lower bound and behaviour of the fundamental solutions; see, e.g., \cite[Section 5.2]{KolV00}. More precisely, our basic hypothesis throughout the paper is the following:
\begin{itemize}
	\item[$(\mathcal{H}_1)$] The spectral measure $\mu$ has a strictly positive density, such that the function $\omega_\mu$ is strictly positive and $(d+1+[\beta])$-times continuously differentiable on $S^{d-1}$.
\end{itemize}
We denote by $(\mathcal{H}_2)$ to refer to $(\mathcal{H}_1)$ whenever we need to assume that $\omega_\mu$ is $(d+2+[\beta])$-times continuously differentiable on $S^{d-1}$, $[\beta]$ being the maximal integer not exceeding the real number $\beta$.

Currently, the theory of pseudo-differential operators is developing because operators like $-\Psi_{\beta}(-i\nabla)$ have an extension generating a Feller semigroup or a sub-Markovian semigroup, under the condition that the symbol $\psi:\RR^d\rightarrow\mathbb{C}$ is a continuous and negative definite function; see \cite[Examples 4.1.12 and 4.1.13]{Jac01}. In our case, $\psi$ satisfies such condition (\cite[Formula 1.9]{Kol00}, \cite[Theorem 3.6.11 and Lemma 3.6.8]{Jac01}). On the other hand, a result of probability theory states that to each Feller semigroup there corresponds a Markov process. Moreover, the connection of the theory of pseudo-differential operators with probability is stronger due to the famous characterization given by the fundamental theorem of Courège (see, e.g., \cite{Kol09,Schi01}). When the stochastic process has an order of stability $\beta$ that is constant, it is called stable jump-diffusions process and the Green function of the Cauchy problem $\dfrac{\partial u}{\partial t}+\Psi_{\beta}(-i\nabla) u=0$ is non-Gaussian, which is interpreted as the transition probability density of the stochastic process \cite[Chapter 7]{Kol11}. Sometimes, these processes are also called Feller processes with pseudo-differential generators, Lévy processes or diffusion with jumps.

The study of \textit{stable non-Gaussian processes} and their generalizations is motivated by the increasing use in the mathematical modeling of processes in engineering, natural sciences and economics. See, e.g., \cite{MJCB14} and \cite[Chapter 1]{Zol86}. An important aspect to mention is the connection between the densities of stable laws
and the \textit{Mittag-Leffler functions} \cite{UZ99}, whose analytical properties are particularly used in this work. Similar to the case of Gaussian processes, which have been widely studied (see, e.g., \cite{Aro67}, \cite{Fri83}, \cite{IKO02}), it arises an interest in the existence, qualitative properties, two-sided estimates and asymptotic behavior for the Green function of non-Gaussian ones (see, e.g., \cite{EK04}, \cite{Kol00}, \cite{KSVZ16}, \cite{KV14}).

In the case $\beta=2$ and $\omega_\mu\equiv 1$ we see that the operator, namely $\Psi_2(-i\nabla)$, becomes the negative Laplacian $(-\Delta)$ with symbol $\psi(\xi)=\norm{\xi}^2$. The corresponding fundamental solution with $\alpha\in (0,1)$ has been studied (see, e.g. \cite[Chapter 5]{EIK04}) and bounds can be found in \cite{EK04}. These bounds are used in \cite{KSVZ16} in order to study the fundamental solution $Z$ for the subdiffusion problem 
\begin{equation*}
\begin{split}
\partial_t^\alpha(u-u_0)-\Delta u &=0,\quad t>0,\; x\in\mathbb{R}^d,\\
u|_{t=0}&=u_0,\quad x\in\mathbb{R}^d.
\end{split}
\end{equation*}
In this work, the authors show the existence of critical values for $Z(t,\cdot)$, $t>0$, concerning $L_p$-integrability on $\mathbb{R}^d$. They also find conditions so that $Z$ belongs to $L_p$-weak space on $\mathbb{R}^d$. For a non-homogeneous analogue Cauchy problem, there exists a \textit{Green matrix} $(Z,Y)$, that is, a pair of fundamental solutions (\cite{EK04}). A more general context for the fractional Laplacian is considered in \cite{KeSZ17} where a representation of $(Z,Y)$ is given by Fox $H$-functions.

Recently in \cite[Section 2]{JK19} and \cite[Section 8.2]{Kol19} the authors show that the linear Cauchy problem
\begin{equation*}
	\begin{split}
		\partial^{\alpha}_t(u-u_0)(t,x)+\Psi_{\beta}(-i\nabla)u(t,x)&=f(t,x),\quad t>0, \; x\in\mathbb{R}^d,\\ u(t,x)|_{t=0}&=u_0(x),\quad x\in\mathbb{R}^d, 
	\end{split}
\end{equation*}
admits a pair of fundamental solutions $(Z,Y)$, given by
\begin{equation}\label{Z:G}
	Z(t,x):=\dfrac{1}{\alpha}\displaystyle\int_0^\infty G(t^\alpha s,x)s^{-1-\frac{1}{\alpha}}G_\alpha(1,s^{-\frac{1}{\alpha}})ds
\end{equation}
and
\begin{equation}\label{Y:G}
	Y(t,x):=\displaystyle\int_0^\infty t^{\alpha-1}G(t^\alpha s,x)s^{-\frac{1}{\alpha}}G_\alpha(1,s^{-\frac{1}{\alpha}})ds,
\end{equation}
where $G$ stands for the Green function that solves the equation
\[
\partial_t\, v(t,x)+\Psi_{\beta}(-i\nabla)v(t,x)=0,\quad t>0, \; x\in\mathbb{R}^d, 
\]
with the initial condition
\[
G(t,x)|_{t=0}=\delta_0(x),\quad x\in\mathbb{R}^d,
\]
$\delta_0$ being the Dirac delta distribution, and $G_\alpha(\cdot,\cdot)$ is the Green function that solves the problem 
\[
\partial_t\, v(t,s)+\dfrac{d^\alpha}{ds^\alpha}v(t,s)=0,\quad t>0,\;s\in\RR, \; G_\alpha(0,s)=\delta(s), 
\]
where $\alpha\in (0,1)$ and 
\[\dfrac{d^\alpha}{ds^\alpha}f(s):=\dfrac{1}{\Gamma(-\alpha)}\int_0^\infty \dfrac{f(s-\tau)-f(s)}{\tau^{1+\alpha}}d\tau,
\] 
see \cite[Formulas (1.111) and (2.74)]{Kol19}. We use the representation of the fundamental solutions $(Z,Y)$ given by \eqref{Z:G}-\eqref{Y:G} to obtain our main results.

Finally, our choice of the non-linear term $\lambda\abs{\cdot}^{\gamma-1}(\cdot)$ in \eqref{general} is a generalization of the classical \textit{rigid ignition model}. The behaviour of the combustion processes, involving non-linear source terms, has become a challenging field for mathematical analysis in the last decades (see, e.g., \cite[Chapter 3]{BE89},\cite{GNZ20}).

This paper is organized as follows. In Sect.~\ref{sec:1} we have compiled some properties of the fundamental solutions. We show that $Z,Y\in C((0,\infty); L_p(\RR^d))$ for all $1\leq p < \kappa_{i}$ ($i=1,2$)  with $\kappa_1$ as in Theorem \ref{CotaNormapZ} and $\kappa_2$ as in Theorem \ref{CotaNormapY} respectively, and they are locally Lipschitz functions in $(0, \infty)$, i.e., for any $\epsilon> 0$ and $s,t\geq\epsilon$ there exists a constant $C>0$ depending on $\epsilon, d,p,\beta, \alpha$, such that
\[
\norm{Z(t,\cdot)-Z(s,\cdot)}_p\leq C |t-s|,
\]
and
\[
\norm{Y(t,\cdot)-Y(s,\cdot)}_p\leq C |t-s|,
\]
respectively. Further, we also prove in this section that $Z(t,\cdot)$ and $\frac{1}{g_{\alpha}(t)}Y(t,\cdot)$ are {\it approximations of the identity} in $L_p$ (for all $1\leq p < \infty$) as $t\to 0$.

Sect.~\ref{sec:2} establishes the existence and uniqueness of a mild solution to \eqref{general}, which is handled by the fixed points theory of operators. We start by presenting a local Banach space of functions and defining a nonlinear operator on it. The main result of this section is stated in Theorem \ref{localsolution} and the detailed proof relies on a delicate analysis of some parameters. In Sect.~\ref{sec:3} we find a global solution using similar arguments as in preceding section, but with a Banach space of functions for large times and a small enough initial condition. For this purpose we have introduced a new parameter, which allows us to use an interpolation result. The main result of this section is given by Theorem \ref{globalsolutionsmallu0}. Sect.~\ref{sec:4} is devoted to set conditions to ensure the existence and uniqueness of a global positive solution, using the theory of Volterra equations and the results from previous sections. We prove a main result in Theorem \ref{globalpositivesolution}, the theory of pseudo-differential operators associated with a negative definite symbol being particularly important for this purpose. In the last section we show that any global solution of \eqref{general} describes an asymptotic profile as $t\rightarrow\infty$ of the form $AZ+BY$ in $L_p$, for some real constants $A$ and $B$, where $p$ and other parameters satisfy suitable conditions as in Theorem \ref{decay}. From our knowledge, such kind of result is unknown in the literature even in the case of $\beta=2$ and $\omega_\mu\equiv 1$. This is due to the fact that $Z$ and $Y$ are different whenever $\alpha\in (0,1)$. From \cite[Theorem 2.7]{HKNS06} it is known that if $\alpha = 1$, $\beta=2$ and $\omega_\mu\equiv 1$, then $Z = Y$ and the corresponding global solution of \eqref{general} behaves as $(A+B)Z$ in $L_p$, as $t\rightarrow\infty$.     

\section{Preliminaries}
\label{sec:1}
In what follows we use the notations $f\asymp g$ and $f\lesssim g$ in $D$, which means that there exists constants  $C,C_1,C_2>0$ such that $C_1 g\leq f\leq C_2 g$ and $f\leq C g$ in $D$, respectively. Such constants may change line by line. We say that $f(x)$ and $g(x)$ are asymptotically equivalent as $x\rightarrow\infty$, if the quotient $\dfrac{f(x)}{g(x)}$ tends to unity. In this case, our notation is $f(x)\sim g(x)$ $(x\rightarrow\infty)$ (see \cite{Bru58}). We also use the notation $\Omega=\norm{x}^\beta t^{-\alpha}$ for $x\in \RR^d$ and $t>0$.

\medbreak

The following result summarize the two-side estimates for $Z$. For its proof see \cite[Theorem 2]{JK19}. 

\begin{proposition}
	\label{cotasZ}
	Let $\alpha\in (0,1)$ and $\beta\in (0,2)$. Assume the hypothesis $(\mathcal{H}_1)$ holds. Then there exists a positive constant $C$ such that for $(t,x)\in(0,\infty)\times\mathbb{R}^d$ the following two-sided estimates for $Z$ hold. For $\Omega\leq 1$,
		\begin{align}
        &Z(t,x)\asymp C t^{-\frac{\alpha d}{\beta}}&\text{if}\quad d<\beta,\label{Z1}\\		
		&Z(t,x)\asymp C t^{-\alpha}(|\log(\Omega)|+1) &\text{if}\quad d=\beta, \label{Z2}\\
		&Z(t,x)\asymp C t^{-\frac{\alpha d}{\beta}}\Omega^{1-\frac{d}{\beta}} &\text{if}\quad d>\beta.\label{Z3}
		\end{align}
		
		For $\Omega\geq 1$,
		\begin{equation}
		\label{Z4}
		Z(t,x)\asymp C t^{-\frac{\alpha d}{\beta}}\Omega^{-1-\frac{d}{\beta}}.
		\end{equation}
\end{proposition}

In the same way we have derived  the two-side estimates for $Y$.

\begin{proposition}
	\label{cotasY}
	Under the same assumptions as Proposition \ref{cotasZ}, the following two-sided estimates for $Y$ hold.
	For $\Omega\leq 1$,
		\begin{align}
		&Y(t,x)\asymp C t^{-\frac{\alpha d}{\beta}+\alpha-1}&\text{if}\quad d<2\beta,\label{Y1}\\
		&Y(t,x)\asymp C t^{-\alpha-1}(|\log(\Omega)|+1)&\text{if}\quad d=2\beta, \label{Y2}\\
		&Y(t,x)\asymp C t^{-\frac{\alpha d}{\beta}+\alpha-1}\Omega^{2-\frac{d}{\beta}}&\text{if}\quad d>2\beta.\label{Y3}
		\end{align}
		
		For $\Omega\geq 1$,
		\begin{equation}
		\label{Y4}
		Y(t,x)\asymp C t^{-\frac{\alpha d}{\beta}+\alpha-1}\Omega^{-1-\frac{d}{\beta}}.
		\end{equation}
\end{proposition}
\begin{proof}
	The assertions follow from a straightforward computations made in the proof of the previous estimates for $Z$.
\end{proof}

\begin{remark}
We note a singularity at the origin with respect to the spatial variable for $Z$, whenever $d\geq\beta$, and for $Y$ whenever $d\geq 2\beta$. It is well known that this type of singularities occurs in the equations of fractional evolution in time, even if $\beta=2$ and $\omega_\mu\equiv 1$.
\end{remark}
Now we continue by showing some properties of the fundamental solutions $Z$ and $Y$.
\begin{lemma}
	\label{cotasDeltaZ}
	Under the same assumptions as Proposition \ref{cotasZ}, there exists a positive constant $C$ for all $t_1,t_2>0$ and $x\in\RR^d$, such that there exists $t_c>0$, between $t_1$ and $t_2$, and the following estimates for $Z$ hold with $\Omega_c=\norm{x}^\beta t_c^{-\alpha}$. For $\Omega_c\leq 1$,		
	\begin{equation*}
		|Z(t_1,x)-Z(t_2,x)|\leq C |t_1-t_2|
		\begin{cases}
			t_c^{-\frac{\alpha d}{\beta}-1} &\quad \text{if}\quad d<\beta,\\
			t_c^{-\alpha-1}(|\log(\Omega_c)|+1) &\quad \text{if}\quad d=\beta,\\
			t_c^{-\frac{\alpha d}{\beta}-1}\Omega_c^{1-\frac{d}{\beta}}&\quad \text{if}\quad d>\beta,		
		\end{cases}
	\end{equation*}
	and for $\Omega_c\geq 1$,
	\[
	|Z(t_1,x)-Z(t_2,x)|\leq C |t_1-t_2| t_c^{-\frac{\alpha d}{\beta}-1}\Omega_c^{-1-\frac{d}{\beta}}.
	\]
\end{lemma}	

\begin{proof}
	From \eqref{Z:G} it follows that
	\[
	Z(t_1,x)-Z(t_2,x)=\dfrac{1}{\alpha}\displaystyle\int_0^\infty \left[G(t_1^\alpha s,x)-G(t_2^\alpha s,x)\right] s^{-1-\frac{1}{\alpha}}G_\alpha(1,s^{-\frac{1}{\alpha}})ds.
	\]
	It is known (\cite[Theorem 4.5.1]{Kol19}) that $G$ is differentiable with respect to $t>0$ and satisfies
	\begin{equation*}
		\begin{split}
			|G(t,x)|\leq C \min\left(t^{-\frac{d}{\beta}},\dfrac{t}{\norm{x}^{d+\beta}}\right),\\
			\left|t\dfrac{\partial G}{\partial t}(t,x)\right |\leq C \min\left(t^{-\frac{d}{\beta}},\dfrac{t}{\norm{x}^{d+\beta}}\right).
		\end{split}
	\end{equation*}
	In these estimates, $C$ depends on $\beta$, $d$ and the bounds for $\omega_\mu$. Using this and the mean-value theorem, we have that for some $t_c$ between $t_1$ and $t_2$,
	\begin{equation*}
		\begin{split}
			&|Z(t_1,x)-Z(t_2,x)|\\
			&\leq \dfrac{1}{\alpha}\displaystyle\int_0^\infty \left|\dfrac{\partial G(t^\alpha s,x)}{\partial t}\right|_{t=t_c}|t_1-t_2| s^{-1-\frac{1}{\alpha}}G_\alpha(1,s^{-\frac{1}{\alpha}})ds\\
			&\leq \dfrac{|t_1-t_2|}{t_c}\displaystyle\int_0^\infty \left|t_c^\alpha s G'(t_c^\alpha s,x)\right| s^{-1-\frac{1}{\alpha}}G_\alpha(1,s^{-\frac{1}{\alpha}})ds\\
			&\leq C \dfrac{|t_1-t_2|}{t_c}\displaystyle\int_0^\infty \min\left((t_c^\alpha s)^{-\frac{d}{\beta}}, \dfrac{t_c^\alpha s}{\norm{x}^{d+\beta}}\right) s^{-1-\frac{1}{\alpha}}G_\alpha(1,s^{-\frac{1}{\alpha}})ds.
		\end{split}
	\end{equation*}
	From \cite[Proposition 2.4.1]{Kol19} and \cite[Theorem 2.5.2]{Zol86} one finds that the asymptotic behavior of $G_\alpha$ is the same as for the density $w_\alpha$ given in \cite[Proposition 1]{JK19} (with the skewness of the distribution that equals to $0$) by
\[
		w_\alpha(\tau) \sim C
		\begin{cases}
			\tau^{-1-\alpha} &\quad \text{as}\quad \tau\rightarrow\infty,\\
			\tau^{-\frac{2-\alpha}{2(1-\alpha)}}e^{-c_\alpha\tau^{-\frac{\alpha}{1-\alpha}}} &\quad \text{as}\quad \tau\rightarrow 0,\end{cases}		
\]	
where $c_\alpha=(1-\alpha)\alpha^{\frac{\alpha}{1-\alpha}}$. Therefore, we proceed in the same way from \cite[Theorem 2]{JK19}, i.e.,
	\begin{align*}
		&\displaystyle\int_0^\infty \min\left((t_c^\alpha s)^{-\frac{d}{\beta}}, \dfrac{t_c^\alpha s}{\norm{x}^{d+\beta}}\right) s^{-1-\frac{1}{\alpha}}G_\alpha(1,s^{-\frac{1}{\alpha}})ds\\
		&=\displaystyle\int_0^1 \min\left((t_c^\alpha s)^{-\frac{d}{\beta}}, \dfrac{t_c^\alpha s}{\norm{x}^{d+\beta}}\right) s^{-1-\frac{1}{\alpha}}w_\alpha(s^{-\frac{1}{\alpha}})ds\\
		&~~~~+\displaystyle\int_1^\infty \min\left((t_c^\alpha s)^{-\frac{d}{\beta}}, \dfrac{t_c^\alpha s}{\norm{x}^{d+\beta}}\right) s^{-1-\frac{1}{\alpha}}w_\alpha(s^{-\frac{1}{\alpha}})ds.
	\end{align*}
With the asymptotic behavior of $w_\alpha$, the first integral reduces to
\[
C\displaystyle\int_0^1 \min\left((t_c^\alpha s)^{-\frac{d}{\beta}}, \dfrac{t_c^\alpha s}{\norm{x}^{d+\beta}}\right) ds
\]  	
and we observe that the improper integral $\displaystyle\int_0^1 (t_c^\alpha s)^{-\frac{d}{\beta}} ds$ appears whenever $\Omega_c <s <1$. Here, we need to check the cases $d=\beta$, $d<\beta$ and $d>\beta$, respectively. Thus, we get the desired bounds.
\end{proof}

\begin{lemma}
	\label{cotasDeltaY}
	Under the same assumptions as Proposition \ref{cotasZ}, there exists a positive constant $C$ for all $t_1,t_2>0$ and $x\in\RR^d$, such that there exists $t_c>0$, between $t_1$ and $t_2$, and the following estimates for $Y$ hold with $\Omega_c=\norm{x}^\beta t_c^{-\alpha}$. For $\Omega_c\leq 1$,
	\[
	|Y(t_1,x)-Y(t_2,x)|\leq C |t_1-t_2|	
	\begin{cases}
		t_c^{-\frac{\alpha d}{\beta}+\alpha-2} &\quad \text{if}\quad d<2\beta,\\
		t_c^{-\alpha-2}(|\log(\Omega_c)|+1) &\quad \text{if}\quad d=2\beta,\\
		t_c^{-\frac{\alpha d}{\beta}+\alpha-2}\Omega_c^{2-\frac{d}{\beta}}&\quad \text{if}\quad d>2\beta,		
	\end{cases}
	\]
	and for $\Omega_c\geq 1$,
	\[
	|Y(t_1,x)-Y(t_2,x)|\leq C |t_1-t_2| t_c^{-\frac{\alpha d}{\beta}+\alpha-2}\Omega_c^{-1-\frac{d}{\beta}}.
	\]	
\end{lemma}
\begin{proof}
	The assertions follow from a straightforward computations made in the proof of the previous estimates for $Z$, but using \eqref{Y:G}.
\end{proof}

\medbreak

\begin{lemma}
	\label{cotasDeltaEspacialZ}
Let $\alpha\in (0,1)$ and $\beta\in (0,2)$. Assume the hypothesis $(\mathcal{H}_2)$ holds. Then there exists a positive constant $C$ for all $t>0$ and $x_1,x_2\in\RR^d$, such that there exists $\zeta$ in the open segment connecting $x_1$ and $x_2$, and the following estimates for $Z$ hold with $\Omega_\zeta=\norm{\zeta}^\beta t^{-\alpha}$. For $\Omega_\zeta\leq 1$,		
    \[
    |Z(t,x_1)-Z(t,x_2)|\leq C \norm{x_1-x_2} t^{-\frac{\alpha(d+1)}{\beta}}\Omega_\zeta^{1-\frac{d+1}{\beta}}
    \]
	and for $\Omega_\zeta\geq 1$,
	\[
	|Z(t,x_1)-Z(t,x_2)|\leq C \norm{x_1-x_2} t^{-\frac{\alpha(d+1)}{\beta}}\Omega_\zeta^{-1-\frac{d+1}{\beta}}.
	\]
\end{lemma}	

\begin{proof}
	From \eqref{Z:G} it follows that
	\[
	Z(t,x_1)-Z(t,x_2)=\dfrac{1}{\alpha}\displaystyle\int_0^\infty \left[G(t^\alpha s,x_1)-G(t^\alpha s,x_2)\right] s^{-1-\frac{1}{\alpha}}G_\alpha(1,s^{-\frac{1}{\alpha}})ds.
	\]
	From (\cite[Theorem 4.5.1]{Kol19}) we know that $G$ is one time continuously differentiable in $x$ and satisfies, for any $j=1,\cdots, d$,
	\begin{equation*}
		\left|\dfrac{\partial G}{\partial x_j}(t,x)\right |\leq C \min\left(t^{-\frac{d+1}{\beta}},\dfrac{t}{\norm{x}^{d+\beta+1}}\right).
		\end{equation*}
	We recall that $C$ depends on $\beta$, $d$ and the bounds for $\omega_\mu$. Let $DG(t^\alpha s,x)$ the Jacobian of $G(t^\alpha s,\cdot)$ in the point $x$. By the mean-value inequality, we have that for some $\zeta$ in the open segment between $x_1$ and $x_2$,
	\begin{equation*}
		\begin{split}
			&|Z(t,x_1)-Z(t,x_2)|\\
			&\leq \frac{1}{\alpha}\displaystyle\int_0^\infty \abs{DG(t^\alpha s,\zeta)(x_1-x_2)} s^{-1-\frac{1}{\alpha}}G_\alpha(1,s^{-\frac{1}{\alpha}})ds\\
			&\leq \frac{\norm{x_1-x_2}}{\alpha}\displaystyle\int_0^\infty \sum_{j=1}^d\left|\dfrac{\partial G}{\partial x_j}(t^\alpha s,\zeta)\right| s^{-1-\frac{1}{\alpha}}G_\alpha(1,s^{-\frac{1}{\alpha}})ds\\
			&\lesssim\norm{x_1-x_2}\displaystyle\int_0^\infty \min\left((t^\alpha s)^{-\frac{d+1}{\beta}}, \dfrac{t^\alpha s}{\norm{\zeta}^{d+\beta+1}}\right) s^{-1-\frac{1}{\alpha}}G_\alpha(1,s^{-\frac{1}{\alpha}})ds.
		\end{split}
	\end{equation*}
	Now, we proceed in the same way from \cite[Theorem 2]{JK19}.
\end{proof}

\begin{lemma}
	\label{cotasDeltaEspacialY}
	Under the same assumptions as Lemma \ref{cotasDeltaEspacialZ}, then there exists a positive constant $C$ for all $t>0$ and $x_1,x_2\in\RR^d$, such that there exists $\zeta$ in the open segment connecting $x_1$ and $x_2$, and the following estimates for $Y$ hold with $\Omega_\zeta=\norm{\zeta}^\beta t^{-\alpha}$. For $\Omega_\zeta\leq 1$,
	\begin{equation*}
		|Y(t,x_1)-Y(t,x_2)|\leq C \norm{x_1-x_2}
		\begin{cases}
			t^{-\frac{\alpha(d+1)}{\beta}+\alpha-1} &\quad \text{if}\quad d+1<2\beta,\\
			t^{-\alpha-1}(|\log(\Omega_\zeta)|+1) &\quad \text{if}\quad d+1=2\beta,\\
			t^{-\frac{\alpha(d+1)}{\beta}+\alpha-1}\Omega_\zeta^{2-\frac{d+1}{\beta}}&\quad \text{if}\quad d+1>2\beta,		
		\end{cases}
	\end{equation*}
	and for $\Omega_\zeta\geq 1$,
	\[
	|Y(t,x_1)-Y(t,x_2)|\leq C \norm{x_1-x_2} t^{-\frac{\alpha(d+1)}{\beta}+\alpha-1}\Omega_\zeta^{-1-\frac{d+1}{\beta}}.
	\]	
\end{lemma}
\begin{proof}
	This is similar to the proof of the previous lemma for $Z$.
\end{proof}

\medbreak

Next, we estimate the $L_p$-norm of $Z$. Let $p\geq 1$ and $t\in(0,\infty)$. We begin by splitting the integral on $\mathbb{R}^d$ according to the conditions for $\Omega$ given by (\ref{Z1})-(\ref{Z4}). That is,
\[
\displaystyle\int_{\mathbb{R}^d}Z^p(t,x)dx=\displaystyle\int_{\{\Omega\geq 1\}}Z^p(t,x)dx + \displaystyle\int_{\{\Omega\leq 1\}}Z^p(t,x)dx.
\]
In the case of $\Omega\geq 1$, the integral on this set has two-sided estimates for all $d\geq 1$ and all $\beta\in(0,2)$, given by (\ref{Z4}). Therefore
\[
\displaystyle\int_{\{\Omega\geq 1\}}Z^p(t,x)dx\asymp \displaystyle\int_{\{\Omega\geq 1\}}t^{-\frac{\alpha d p}{\beta}}\Omega^{(-1-\frac{d}{\beta})p} dx.
\]
Setting $r=\norm{x}$, we obtain
\begin{equation*}
\begin{split}
\displaystyle\int_{\{\Omega\geq 1\}}Z^p(t,x)dx&\asymp \displaystyle\int_{t^{\frac{\alpha}{\beta}}}^\infty t^{-\frac{\alpha d p}{\beta}}(r^\beta t^{-\alpha})^{(-1-\frac{d}{\beta})p}r^{d-1} dr\\
&= \displaystyle\int_{t^{\frac{\alpha}{\beta}}}^\infty t^{\alpha p}r^{(-\beta-d)p+d-1} dr\\
&= \displaystyle\int_{t^{\frac{\alpha}{\beta}}}^\infty t^{\alpha p}(t^{\frac{\alpha}{\beta}}\underset{s}{\underbrace{t^{-\frac{\alpha}{\beta}} r}})^{(-\beta-d)p+d-1} dr\\
&= \displaystyle\int_1^\infty t^{-\frac{\alpha dp}{\beta}+\frac{\alpha d}{\beta}}s^{-\beta p-1-(p-1)d} ds\\
&= t^{-\frac{\alpha d p}{\beta}\left(1-\frac{1}{p}\right)} \displaystyle\int_1^\infty s^{-\beta p-1-(p-1)d} ds.
\end{split}
\end{equation*}

The last integral converges if and only if $\beta p+1+(p-1)d>1$, which holds true for all $d\geq 1$, $\beta\in (0,2)$ and $1\leq p <\infty$. Consequently, we obtain the estimate 
\begin{equation}
\label{cotageneralZ}
\displaystyle\int_{\{\Omega\geq 1\}}Z^p(t,x)dx\asymp t^{-\frac{\alpha d p}{\beta}\left(1-\frac{1}{p}\right)}.
\end{equation}

\medbreak

Now, for $\Omega\leq 1$, we consider the following cases separately: $\beta\in (0,1)$, $\beta=1$ and $\beta\in (1,2)$. For $\beta\in (0,1)$, we employ the bounds given by (\ref{Z3}) and we set again $r=\norm{x}$ and the substitution $r=t^{\frac{\alpha}{\beta}} s$, obtaining
\begin{equation*}
\begin{split}
\displaystyle\int_{\{\Omega\leq 1\}}Z^p(t,x)dx&\asymp \displaystyle\int_{\{\Omega\leq 1\}}t^{-\frac{\alpha d p}{\beta}}\Omega^{(1-\frac{d}{\beta})p} dx\\
&\asymp \displaystyle\int_{0}^{t^{\frac{\alpha}{\beta}}} t^{-\frac{\alpha d p}{\beta}}(r^\beta t^{-\alpha})^{(1-\frac{d}{\beta})p}r^{d-1} dr\\
&= \displaystyle\int_{0}^{t^{\frac{\alpha}{\beta}}} t^{-\alpha p}r^{(\beta-d)p+d-1} dr\\
&= \displaystyle\int_0^1 t^{-\frac{\alpha dp}{\beta}+\frac{\alpha d}{\beta}}s^{(\beta-d)p+d-1} ds\\
&= t^{-\frac{\alpha dp}{\beta}\left(1-\frac{1}{p}\right)} \displaystyle\int_0^1 s^{(\beta-d)p+d-1} ds.
\end{split}
\end{equation*}
The last integral converges if and only if $1-d+(d-\beta)p<1$,  which is equivalent with $p<\dfrac{d}{d-\beta}$. This together with (\ref{cotageneralZ}) gives
\begin{equation}
\label{normaZcaso1}
\norm{Z(t,\cdot)}_p\asymp t^{-\frac{\alpha d}{\beta}\left(1-\frac{1}{p}\right)},~ 1\leq p<\dfrac{d}{d-\beta}, ~ d\geq 1,~ \beta\in (0,1).
\end{equation}

For $\beta=1$, we check out $d=\beta=1$ and $d>\beta$. We employ the bounds given by (\ref{Z2}) and (\ref{Z3}), respectively. In the case of $d=1$, we obtain
\begin{equation*}
\begin{split}
\displaystyle\int_{\{\Omega\leq 1\}}Z^p(t,x)dx&\asymp \displaystyle\int_{\{\Omega\leq 1\}}t^{-\alpha p}\left(|\log\Omega|+1\right)^p dx\\
&\asymp \displaystyle\int_{0}^{t^{\alpha}} t^{-\alpha p}\left(|\log(rt^{-\alpha})|+1\right)^p dr\\
&=\displaystyle\int_{0}^{1} t^{-\alpha p+\alpha}\left(|\log(s)|+1\right)^p ds\\
&=t^{-\alpha p\left(1-\frac{1}{p}\right)}\displaystyle\int_{0}^{1} \left(|\log(s)|+1\right)^p ds.
\end{split}
\end{equation*}

The last integral converges for all $1\leq p <\infty$. Together with (\ref{cotageneralZ}) yields
\begin{equation}
\label{normaZcaso2_1}
\norm{Z(t,\cdot)}_p\asymp t^{-\alpha\left(1-\frac{1}{p}\right)},~ 1\leq p<\infty, ~ d=1,~ \beta=1.
\end{equation}

Now, for $d\geq 2$ the corresponding estimate is similar to the case of $\beta\in (0,1)$. Thus, 
\begin{equation}
\label{normaZcaso2_2}
\norm{Z(t,\cdot)}_p\asymp t^{-\alpha d\left(1-\frac{1}{p}\right)},~ 1\leq p<\dfrac{d}{d-1}, ~ d\geq 2,~ \beta=1.
\end{equation}

Finally, for $\beta\in (1,2)$, we check out $d<\beta$ and $d >\beta$. We employ the bounds given by (\ref{Z1}) and (\ref{Z3}), respectively. For $d=1$,
\begin{equation*}
\begin{split}
\displaystyle\int_{\{\Omega\leq 1\}}Z^p(t,x)dx&\asymp \displaystyle\int_{\{\Omega\leq 1\}}t^{-\frac{\alpha d p}{\beta}}dx\\
&\asymp \displaystyle\int_{0}^{t^{\frac{\alpha}{\beta}}} t^{-\frac{\alpha p}{\beta}} dr\\
&= t^{-\frac{\alpha p}{\beta}} t^{\frac{\alpha}{\beta}}\\
&= t^{-\frac{\alpha p}{\beta}\left(1-\frac{1}{p}\right)}.
\end{split}
\end{equation*}

Together with (\ref{cotageneralZ}) and the fact that
\[
\displaystyle\sup_{x\in\mathbb{R}} Z(t,\cdot)\asymp t^{-\frac{\alpha}{\beta}},
\]
we conclude
\begin{equation}
\label{normaZcaso3_1}
\norm{Z(t,\cdot)}_p\asymp t^{-\frac{\alpha}{\beta}\left(1-\frac{1}{p}\right)},~ 1\leq p\leq\infty, ~ d=1,~ \beta\in (1,2).
\end{equation}

In the case of $d\geq 2$, the corresponding estimate is similar to that of $\beta\in (0,1)$. Thus,  
\begin{equation}
\label{normaZcaso3_2}
\norm{Z(t,\cdot)}_p\asymp t^{-\frac{\alpha d}{\beta}\left(1-\frac{1}{p}\right)},~ 1\leq p<\dfrac{d}{d-\beta}, ~ d\geq 2,~ \beta\in (1,2).
\end{equation}

Gathering the two-side estimates from (\ref{normaZcaso1}) to (\ref{normaZcaso3_2}), we have proved the following result.

\begin{theorem}
\label{CotaNormapZ}
Let $d\in\mathbb{N}$, $\alpha\in (0,1)$ and $\beta\in (0,2)$. Assume the hypothesis $(\mathcal{H}_1)$ holds. The kernel $Z(t,\cdot)$ belongs to $L_p(\mathbb{R}^d)$ for all $t>0$ if, and only if, $1\leq p<\kappa_1$, where
\[
\kappa_1=\kappa_1(d,\beta):=
\begin{cases}
\frac{d}{d-\beta} & \text{if } d>\beta,\\
\infty & \text{otherwise}.
\end{cases}
\]
Moreover, the two-side estimate 
\begin{equation}
\label{cotasZp}
\norm{Z(t,\cdot)}_p\asymp t^{-\frac{\alpha d}{\beta}\left(1-\frac{1}{p}\right)},\; t>0
\end{equation}
holds for every $1\leq p<\kappa_1$. In the case of $d<\beta$, (\ref{cotasZp}) remains true for $p=\infty$.
\end{theorem}

We next examine the critical case $p=\frac{d}{d-\beta}$ for $d>\beta$, in the $L_{p}$-weak space with the quasi-norm $|\cdot|_{p,\infty}$ defined by
\[
|f|_{p,\infty}:=\displaystyle\sup_{\lambda>0}\{\lambda~ d_{f}(\lambda)^{\frac{1}{p}}\},
\]
where
\[
d_{f}(\lambda)=|\{x\in\mathbb{R}^d:f(x)>\lambda\}|
\]
stands for the distribution function of $f$.

\begin{theorem} Let $d\in\mathbb{N}$, $\alpha\in (0,1)$ and $\beta\in (0,2)$ such that $d>\beta$. Assume the hypothesis $(\mathcal{H}_1)$ holds. Then $Z(t,\cdot)\in L_{\frac{d}{d-\beta},\infty}(\mathbb{R}^d)$ and satisfies
\[
|Z(t)|_{\frac{d}{d-\beta},\infty}\lesssim t^{-\alpha}, \; t>0. 
\]
\end{theorem}

\begin{proof}
Let $t>0$ and denote $Z(t)=Z(t,\cdot)$. Set $p=\dfrac{d}{d-\beta}$. By definition, $Z(t)\in L_{p,\infty}(\RR^d)$ if the quasi-norm
\[
|Z(t)|_{p,\infty}=\displaystyle\sup_{\lambda>0}\{\lambda~ d_{Z(t)}(\lambda)^{\frac{1}{p}}\}<\infty.
\]
As above we use the similarity variable $\Omega=\norm{x}^\beta t^{-\alpha}$ and we split $Z(t)$ as $Z(t)=Z(t)\chi(t)_{\{\Omega\leq 1\}}+Z(t)\chi(t)_{\{\Omega\geq 1\}}$. Then 
\[
|Z(t)|_{p,\infty}\leq 2\left(|Z(t)\chi(t)_{\{\Omega\leq 1\}}|_{p,\infty}+|Z(t)\chi(t)_{\{\Omega\geq 1\}}|_{p,\infty}\right).
\]
By (\ref{cotageneralZ}) and the $L_p$ version Tchebyshev's inequality (\cite[Formula (5.49)]{WZ77}), we obtain
\begin{equation*}
|Z(t)\chi(t)_{\{\Omega\geq 1\}}|_{p,\infty}\leq \norm{Z(t)\chi(t)_{\{\Omega\geq 1\}}}_{p}\leq C t^{-\frac{\alpha d}{\beta}\left(1-\frac{1}{p}\right)}=Ct^{-\alpha}.
\end{equation*}
Now, employing (\ref{Z3}) we have 
\begin{equation*}
\begin{split}
d_{Z(t)\chi(t)_{\{\Omega\leq 1\}}}(\lambda)&=|\{x\in\mathbb{R}^d:Z(t,x)>\lambda\quad\text{y}\quad\Omega\leq 1\}|\\
&\leq|\{x\in\mathbb{R}^d:\lambda < Ct^{-\frac{\alpha d}{\beta}}\Omega^{1-\frac{d}{\beta}} \}|\\
&=|\{x\in\mathbb{R}^d:\lambda < Ct^{-\frac{\alpha d}{\beta}}(\norm{x}^\beta t^{-\alpha})^{1-\frac{d}{\beta}} \}|\\
&=|\{x\in\mathbb{R}^d:\lambda < Ct^{-\alpha}\norm{x}^{\beta-d} \}|\\
&\leq C_1\left(t^{-\alpha}\lambda^{-1}\right)^{\frac{d}{d-\beta}}.\\
\end{split}
\end{equation*}
Thereby, we find that
\[
\lambda\; d_{Z(t)\chi(t)_{\{\Omega\leq 1\}}}(\lambda)^{\frac{1}{p}}\leq C_1^{\frac{1}{p}} t^{-\alpha},
\]
and thus 
\begin{equation*}
|Z(t)\chi(t)_{\{\Omega\leq 1\}}|_{p,\infty}\lesssim  t^{-\alpha}.
\end{equation*}
The proof is complete.
\end{proof}

Analogously, to the analysis done for the $L_p$-integrability of $Z$, we can obtain the corresponding results for $Y$ using the bounds given by \eqref{Y1}-\eqref{Y4}. 
\begin{theorem}
	\label{CotaNormapY}
	Let $d\in\mathbb{N}$, $\alpha\in (0,1)$ and $\beta\in (0,2)$. Assume the hypothesis $(\mathcal{H}_1)$ holds. The kernel $Y(t,\cdot)$ belongs to $L_p(\mathbb{R}^d)$ for all $t>0$ if, and only if, $1\leq p<\kappa_2$, where
	\[
	\kappa_2=\kappa_2(d,\beta):=
	\begin{cases}
		\frac{d}{d-2\beta} & \text{if } d>2\beta,\\
		\infty & \text{otherwise}.
	\end{cases}
	\]
	Moreover, the two-side estimate 
	\begin{equation}
		\label{cotasYp}
		\norm{Y(t,\cdot)}_p\asymp t^{-\frac{\alpha d}{\beta}\left(1-\frac{1}{p}\right)+(\alpha-1)},\; t>0
	\end{equation}
	holds for every $1\leq p<\kappa_2$. In the case of $d<2\beta$, (\ref{cotasYp}) remains true for $p=\infty$.
\end{theorem}

\begin{theorem} Let $d\in\mathbb{N}$, $\alpha\in (0,1)$ and $\beta\in (0,2)$ such that $d>2\beta$. Assume the hypothesis $(\mathcal{H}_1)$ holds. Then $Y(t,\cdot)\in L_{\frac{d}{d-2\beta},\infty}(\mathbb{R}^d)$ and satisfies
	\[
	|Y(t)|_{\frac{d}{d-2\beta},\infty}\lesssim t^{-\alpha-1}, \; t>0. 
	\]
\end{theorem}

Now, we can establish that $Z(t,\cdot)$ is a locally Lipschitz function in $L_p$ and it is an approximation of the identity. To this end, we first show that $Z$ and $Y$ satisfy the following scaling property.
\begin{lemma}
\label{scalingZeY}
Let $t>0$, $x\in\RR^d$. Then
\[
Z(t,x)  = t^{-\frac{\alpha d}{\beta}}Z\left(1,t^{-\frac{\alpha}{\beta}}x\right)
\]
and
\[
Y(t,x)=t^{-\frac{\alpha d}{\beta}+\alpha-1}Y\left(1,t^{-\frac{\alpha}{\beta}}x\right).
\]
\end{lemma}

\begin{proof}
It is well-known that the Fourier transform in the spatial variable of $Z$, denoted by $\hat{Z}(t,\xi)$, is
\begin{equation}
\label{FourierZ}
\hat{Z}(t,\xi)=E_{\alpha,1}(-t^{\alpha}\psi(\xi)),\quad t>0,
\end{equation}
where $E_{\alpha,\delta}(z) :=\displaystyle\sum_{k=0}^\infty\dfrac{z^k}{\Gamma(k\alpha+\delta)}$ stands for the Mittag-Leffler function with two parameters $\alpha, \delta > 0$;  $z\in\mathbb{C}$ (see, e.g., \cite[Sub-section 8.4]{KeSZ17}, \cite[Section 5]{MG00}). In this case the series converges for all values of $z$, so $E_{\alpha,\delta}$ is an entire function. Using this and changing the integration variables, where $\overline{\xi}=\frac{\xi}{\norm{\xi}}$ and $\omega_\mu$ is given by \eqref{omega_mu}, we get
\begin{align*}
&Z(t,x)\\
&=\dfrac{1}{(2\pi)^d}\displaystyle\int_{\mathbb{R}^d}e^{i\, x\cdot\xi} E_{\alpha,1}(-t^{\alpha}\psi(\xi))d\xi\\
&=\dfrac{1}{(2\pi)^d}\displaystyle\int_0^\infty\int_{S^{d-1}}e^{i\left(\dfrac{x}{t^{\frac{\alpha}{\beta}}}\cdot\overline{\xi}\right)t^{\frac{\alpha}{\beta}}\norm{\xi}} E_{\alpha,1}(-(t^{\frac{\alpha}{\beta}}\norm{\xi})^{\beta}\omega_\mu(\overline{\xi}))\norm{\xi}^{d-1}\theta(d\overline{\xi})d\norm{\xi}\\
&=\dfrac{1}{(2\pi)^d}\displaystyle\int_0^\infty\int_{S^{d-1}}e^{i\left(\dfrac{x}{t^{\frac{\alpha}{\beta}}}\cdot\overline{\xi}\right)t^{\frac{\alpha}{\beta}}\norm{\xi}} E_{\alpha,1}(-(t^{\frac{\alpha}{\beta}}\norm{\xi})^{\beta}\omega_\mu(\overline{\xi}))(t^{-\frac{\alpha}{\beta}}\underset{r}{\underbrace{t^{\frac{\alpha}{\beta}}\norm{\xi}}})^{d-1}\theta(d\overline{\xi})d\norm{\xi}\\
&=t^{-\frac{\alpha d}{\beta}}\dfrac{1}{(2\pi)^d}\displaystyle\int_0^\infty\int_{S^{d-1}}e^{i\left(\dfrac{x}{t^{\frac{\alpha}{\beta}}}\cdot\overline{\xi}\right)r} E_{\alpha,1}(-r^{\beta}\omega_\mu(\overline{\xi}))r^{d-1}\theta(d\overline{\xi})dr\\
&=t^{-\frac{\alpha d}{\beta}}Z\left(1,t^{-\frac{\alpha}{\beta}}x\right).
\end{align*}
On the other hand, it is also known that the Fourier transform in the spatial variable of $Y$ is given by
\begin{equation}
\label{FourierY}
\hat{Y}(t,\xi)=t^{\alpha-1}E_{\alpha,\alpha}(-t^{\alpha}\psi(\xi)),\quad t>0.
\end{equation}
Therefore, the previous argument can be applied to 
\[
Y(t,x)=\dfrac{1}{(2\pi)^d}\displaystyle\int_{\mathbb{R}^d}e^{i\, x\cdot\xi} t^{\alpha-1}E_{\alpha,\alpha}(-t^{\alpha}\psi(\xi))d\xi.
\]
\end{proof}

\begin{theorem}
\label{C1menosZ}
Let $d\in\mathbb{N}$, $\alpha\in (0,1)$, $\beta\in (0,2)$. Assume the hypothesis $(\mathcal{H}_1)$ holds. Then
\begin{enumerate}
\item $Z\in C((0,\infty); L_p(\RR^d))$ for $1\leq p<\kappa_1$. Further, for each $\epsilon>0$ there exists a constant $C>0$ depending on $\epsilon, d, p,\alpha,\beta$ such that
\begin{equation}
\label{cotasC1menosZp}
\norm{Z(t,\cdot)-Z(s,\cdot)}_p\leq C|t-s|,
\end{equation}
holds for all $t,s\geq \epsilon$. In the case of $d<\beta$, (\ref{cotasC1menosZp}) remains true for $p=\infty$.

\item For any $v\in L_p(\RR^d)$, with $1\leq p <\infty$, we have 
\[
\lim_{t\to 0}\norm{Z(t,\cdot) \star v - v}_p = 0. 
\]
\end{enumerate}
\end{theorem}

\begin{proof}
Let $0<\epsilon\leq s,t$. From Lemma \ref{cotasDeltaZ} we know that there exists $\tau>0$ between $s$ and $t$, with $\Omega_c=\norm{x}^\beta \tau^{-\alpha}$, such that 
	\begin{equation*}
	|Z(t,x)-Z(s,x)|\leq C |t-s|
	\begin{cases}
	\tau^{-\frac{\alpha d}{\beta}-1} &\quad \text{if}\quad d<\beta,\\
	\tau^{-\alpha-1}(|\log(\Omega_c)|+1) &\quad \text{if}\quad d=\beta,\\
	\tau^{-\frac{\alpha d}{\beta}-1}\Omega_c^{1-\frac{d}{\beta}}&\quad \text{if}\quad d>\beta,		
	\end{cases}
	\end{equation*}
for $\Omega_c\leq 1$ and
	\[
	|Z(t,x)-Z(s,x)|\leq C |t-s| \tau^{-\frac{\alpha d}{\beta}-1}\Omega_c^{-1-\frac{d}{\beta}}
	\]
for $\Omega_c\geq 1$. 
Similar arguments as in the proof of Theorem \ref{CotaNormapZ}, with the only difference that there appears $\Omega$ instead of $\Omega_c$ in the spatial integral respect to $x$, show that
\[
\norm{Z(t,\cdot)-Z(s,\cdot)}_p\leq C_1 |t-s| \tau^{-\frac{\alpha d}{\beta}\left(1-\frac{1}{p}\right)-1}
\]
if and only if $1\leq p<\kappa_1$. Recall that the condition \textit{if and only if} guarantees the existence of the improper Riemann integrals in the proof.

In order to get \eqref{cotasC1menosZp}, we use the fact that $\tau>\epsilon$. Thus,
\[
\norm{Z(t,\cdot)-Z(s,\cdot)}_p\leq C_1 |t-s| \epsilon^{-\frac{\alpha d}{\beta}\left(1-\frac{1}{p}\right)-1}
\]
for all $t,s\geq \epsilon$ and we take $C=C_1\epsilon^{-\frac{\alpha d}{\beta}\left(1-\frac{1}{p}\right)-1}$. This proves $(i)$.

To prove $(ii)$, let $v\in L_p(\RR^d)$ with $p\geq 1$. From \eqref{FourierZ} and since $\psi(0)=0$, we have
\begin{equation}\label{IntegralZ:1}
\int_{\RR^d}Z(t,x)dx = \hat{Z}(t,0) = 1,\quad t>0.
\end{equation}

Now, for $t>0$ define $\phi_t(x):=t^{-\frac{\alpha d}{\beta}}Z\left(1,t^{-\frac{\alpha}{\beta}}x\right)$, $x\in\RR^d$. From \eqref{IntegralZ:1} it follows that $\phi_t\star v\in L_p (\RR^d)$ for all $p\geq 1$.

By applying the Minkowski's integral inequality we obtain 
\begin{equation*}
\begin{split}
\norm{\phi_t\star v-v}_p&=\left(\int_{\RR^d}|(\phi_t\star v)(x)-v(x)|^p dx\right)^{\frac{1}{p}}\\
&=\left(\int_{\RR^d}\left|\int_{\RR^d}v(x-y)\phi_t(y)dy-\int_{\RR^d}v(x)\phi_t(y)dy\right|^p dx\right)^{\frac{1}{p}}\\
&=\left(\int_{\RR^d}\left|\int_{\RR^d}(v(x-y)-v(x))\phi_t(y)dy\right|^p dx\right)^{\frac{1}{p}}\\
&=\left(\int_{\RR^d}\left|\int_{\RR^d}(v(x-t^{\frac{\alpha}{\beta}}y)-v(x))Z(1,y)dy\right|^p dx\right)^{\frac{1}{p}}\\
&\leq\int_{\RR^d}\left(\int_{\RR^d}|v(x-t^{\frac{\alpha}{\beta}}y)-v(x)|^p Z^p(1,y) dx \right)^{\frac{1}{p}}dy\\
&=\int_{\RR^d}Z(1,y)\norm{v(\cdot-t^{\frac{\alpha}{\beta}}y)-v(\cdot)}_p\;dy.
\end{split}
\end{equation*}
Since $\norm{v(\cdot-t^{\frac{\alpha}{\beta}}y)-v(\cdot)}_p \to 0$ as $t\to 0$ and $\norm{v(\cdot-t^{\frac{\alpha}{\beta}}y)-v(\cdot)}_p \leq 2\norm{v}_p$, we apply the dominated convergence theorem to the last integral concluding $Z(t,\cdot)\star v \to v$ in $L_p(\RR^d)$ as $t\to 0$.
\end{proof}

In the same way, we can establish that $Y(t,\cdot)$ is locally Lipschitz function in $L_p$ and that $\frac{1}{g_{\alpha}(t)}Y(t,\cdot)$ is an approximation of the identity in $L_p$ as $t\to 0$.

\begin{theorem}
\label{C1menosY}
Let $d\in\mathbb{N}$, $\alpha\in (0,1)$, $\beta\in (0,2)$. Assume the hypothesis $(\mathcal{H}_1)$ holds. Then
\begin{enumerate}
\item $Y\in C((0,\infty); L_p(\RR^d))$ for $1\leq p<\kappa_2$. Further, for each $\epsilon>0$ there exists a constant $C>0$ depending on $\epsilon, d, p,\alpha,\beta$ such that
\begin{equation}
\label{cotasC1menosYp}
\norm{Y(t,\cdot)-Y(s,\cdot)}_p\leq C|t-s|,
\end{equation}
holds for all $t,s\geq \epsilon$. In the case of $d<2\beta$, (\ref{cotasC1menosYp}) remains true for $p=\infty$.

\item For any $v\in L_p(\RR^d)$, with $1\leq p < \infty$, we have that
\[
\lim_{t\to 0}\left\|\frac{1}{g_{\alpha}(t)}Y(t,\cdot) \star v - v\right\|_p = 0. 
\]
\end{enumerate}
\end{theorem}
\begin{proof}
The proof of $(i)$ is similar to the one we used in Theorem \ref{C1menosZ} part $(i)$. The key point to prove $(ii)$ is Lemma \ref{scalingZeY} and \eqref{FourierY}, because 
\begin{equation*}
	\begin{split}
		\frac{1}{g_\alpha(t)}\displaystyle\int_{\RR^d}Y(t,x)dx=1,\quad t>0.
	\end{split}
\end{equation*}
\end{proof}

We finish this section by proving the following result.
\begin{lemma}
\label{YenfunciondeZ}
Let fixed $x\in\RR^d\setminus \{0\}$. Under the same assumptions as Proposition \ref{cotasZ}, $Z$ and $Y$ satisfy
\[
Y(\cdot,x)=\dfrac{d}{dt}(g_{\alpha}*Z(\cdot,x)),\quad t>0.
\]
\end{lemma}
\begin{proof}
From \eqref{FourierY} and Fubini, we have that
\begin{equation*}
\begin{split}
(g_{1-\alpha}\,*\,Y(\cdot,x))(t)&=\displaystyle\int_0^t g_{1-\alpha}(t-s)Y(s,x)ds\\
&=\displaystyle\int_0^t \dfrac{(t-s)^{-\alpha}}{\Gamma(1-\alpha)}\left[\dfrac{1}{(2\pi)^d}\displaystyle\int_{\mathbb{R}^d}e^{ix\cdot \xi}s^{\alpha-1}E_{\alpha,\alpha}(-\psi(\xi)s^\alpha)d\xi \right] ds\\
&=\dfrac{1}{\Gamma(1-\alpha)(2\pi)^d}\displaystyle\int_{\mathbb{R}^d}e^{ix\cdot \xi}\displaystyle\int_0^t \dfrac{s^{\alpha-1}}{(t-s)^\alpha} E_{\alpha,\alpha}(-\psi(\xi)s^\alpha)  ds\, d\xi.
\end{split}
\end{equation*}
On the other hand, using the definition of $E_{\alpha,\alpha}$ it follows that
\[
\dfrac{s^{\alpha-1}}{(t-s)^\alpha}  E_{\alpha,\alpha}(-\psi(\xi)s^\alpha)=\displaystyle\sum_{k=0}^\infty\dfrac{(-\psi(\xi))^k s^{k\alpha+\alpha-1}}{\Gamma(k\alpha+\alpha)}\dfrac{1}{(t-s)^\alpha}.
\]
By integrating respect to $s$, we obtain
\[
\displaystyle\int_0^t\dfrac{s^{\alpha-1}}{(t-s)^\alpha}  E_{\alpha,\alpha}(-\psi(\xi)s^\alpha)ds=\displaystyle\sum_{k=0}^\infty\dfrac{(-\psi(\xi))^k }{\Gamma(k\alpha+\alpha)}\displaystyle\int_0^t\dfrac{s^{k\alpha+\alpha-1}}{(t-s)^\alpha}ds.
\]
In the last integral, the substitution $s=t\tau$ yields
\[
\displaystyle\int_0^t\dfrac{s^{k\alpha+\alpha-1}}{(t-s)^\alpha}ds=\displaystyle\int_0^1\dfrac{t^{k\alpha}\tau^{k\alpha+\alpha-1}}{(1-\tau)^\alpha}d\tau.
\]
We note that this improper integral exists because $\alpha\in (0,1)$. By using Beta $B$ and Gamma $\Gamma$ functions, we get
\[
\displaystyle\int_0^1 \tau^{k\alpha+\alpha-1}(1-\tau)^{1-\alpha-1}d\tau=B(k\alpha+\alpha,1-\alpha)=\dfrac{\Gamma(k\alpha+\alpha)\Gamma(1-\alpha)}{\Gamma(k\alpha+1)}
\]
and thus
\begin{align*}
\displaystyle\int_0^t\frac{s^{\alpha-1}}{(t-s)^\alpha} E_{\alpha,\alpha}(-\psi(\xi)s^\alpha)ds&=\displaystyle\sum_{k=0}^\infty\dfrac{(-\psi(\xi))^k }{\Gamma(k\alpha+\alpha)}t^{k\alpha}\frac{\Gamma(k\alpha+\alpha)\Gamma(1-\alpha)}{\Gamma(k\alpha+1)}\\
&=\Gamma(1-\alpha)E_{\alpha,1}(-\psi(\xi)t^\alpha).
\end{align*}
Using this in the first part of the proof, we conclude that
\[
Z(t,x)=(g_{1-\alpha}*Y(\cdot,x))(t).
\]
The convolution with $g_\alpha$ and the derivative w.r.t. the time complete the proof.
\end{proof}

\section{Local well-posedness}
\label{sec:2}
Due to the properties of the pair $(g_\alpha,g_{1-\alpha})$, with $0<\alpha<1$ (see, e.g., \cite[Section 1]{KSVZ16}), and under suitable conditions, the problem \eqref{general} can be rewritten as the \textit{semilinear Volterra equation} 
\[
u+g_\alpha*\Psi_{\beta}(-i\nabla)u=u_0+\lambda~g_\alpha * \abs{u}^{\gamma-1}u.
\]
This fact is particularly exploited in Section \ref{sec:4}. However, in this section we deal with the corresponding integral representation for \textit{mild solutions in the sense of Volterra}, which leads to fixed points of the integral equation
	\begin{equation}\label{Inteq:1}
	u(t,x) = \int_{\RR^d} Z(t,x-y)u_0(y)dy + \lambda\int_0^t\int_{\RR^d} Y(t-s,x-y)|u(s,y)|^{\gamma-1}u(s,y)dyds.
	\end{equation}
	By a local solution of the Cauchy problem \eqref{general} we understand the solution $u$ of the corresponding integral equation \eqref{Inteq:1} (the so-called {\textit{mild solution}}) belonging to the Banach space
	\[
	E_T:=C([0,T]; L_p(\RR^d) \cap L_1(\RR^d)) \cap C((0,T]; L_{\infty}(\RR^d)),
	\]
	with the norm
	\[
	\left\|v\right\|_{E_T}:=\sup_{t\in[0,T]}\left(\left\|v(t,\cdot)\right\|_p + \left\|v(t,\cdot)\right\|_1\right) + \sup_{t\in(0,T]}t^{\frac{\alpha d}{\beta p}}\left\|v(t,\cdot)\right\|_{\infty}.
	\]
We define on $E_T$ the operator $\mathcal{M}$ given by 
\begin{equation}
\label{operator}
\mathcal{M}(v)(t,x):=\int_{\RR^d} Z(t,x-y)u_0(y)dy + \lambda\int_0^t\int_{\RR^d} Y(t-s,x-y)|v(s,y)|^{\gamma-1}v(s,y)dyds
\end{equation}
where $u_0\in L_p(\RR^d) \cap L_1(\RR^d)$ is a given data and $v\in E_T$. The space $L_p(\RR^d) \cap L_1(\RR^d)$ is equipped with the usual norm $\norm{\cdot}_1+\norm{\cdot}_p$.

We also need to define the number $\kappa:=\begin{cases}\frac{d}{\beta}, & d>\beta,\\ 1, & otherwise. \end{cases}$

\begin{theorem}
\label{localsolution}
Let $\alpha\in (0,1)$ and $\beta\in (0,2)$. Assume the hypothesis $(\mathcal{H}_1)$ holds. Let $\lambda\in \RR$ and $\gamma>1$. Suppose that $\max\left(1,\kappa,\frac{d(\gamma-1)}{\beta}\right)<p<\infty$. If $u_0\in L_p(\RR^d) \cap L_1(\RR^d)$, then for some $0<T^*<T$ the operator $\mathcal{M}$ defined by \eqref{operator} has a unique fixed point in $E_{T^*}$.  
\end{theorem}

\begin{proof}
Since $u_0\in L_p(\RR^d)\cap L_1(\RR^d)$ and $Z(t,\cdot)\in L_1(\RR^d)$, it follows that $Z(t,\cdot)\star u_0\in L_p(\RR^d)\cap L_1(\RR^d)$ for each $t\in [0,T]$. Condition $\max(1,\kappa)<p<\infty$ and Young’s inequality for convolutions imply that $Z(t,\cdot)\star u_0\in L_{\infty}(\RR^d)$ for each $t\in (0,T]$, since there exists a $q$ such that $\frac{1}{p}+\frac{1}{q}=1$ and $1 < q<\kappa_1$, thus $Z(t,\cdot)\in L_q(\RR^d)$ (see Theorem \ref{CotaNormapZ}). Furthermore,
\begin{equation}
\label{cotasiniciales}
\begin{split}
&\left\|Z(t,\cdot)\star u_0\right\|_1\leq\left\|Z(t,\cdot)\right\|_1\left\|u_0\right\|_1=\left\|u_0\right\|_1,\\
&\left\|Z(t,\cdot)\star u_0\right\|_p\leq\left\|Z(t,\cdot)\right\|_1\left\|u_0\right\|_p=\left\|u_0\right\|_p,\\
&\left\|Z(t,\cdot)\star u_0\right\|_{\infty}\leq\left\|Z(t,\cdot)\right\|_q\left\|u_0\right\|_{p}\lesssim t^{-\frac{\alpha d}{\beta p}}\left\|u_0\right\|_p.
\end{split}
\end{equation}
The continuity of $t\mapsto Z(t,\cdot)\star u_0$ in $[0,T]$ with respect to the norm topology on $ L_p(\RR^d) \cap L_1(\RR^d)$ follows from Theorem \ref{C1menosZ}. The continuity in $(0,T]$ with respect to the norm topology on $L_{\infty}(\RR^d)$ follows from the same considerations mentioned above for $p$ and from Theorem \ref{C1menosZ} part $(i)$.

Let us consider $t\mapsto \displaystyle\int_0^tY(t-s,\cdot)\star|v(s,\cdot)|^{\gamma-1}v(s,\cdot)ds$, where $0<t\leq T$.

Previous arguments, together with the Minkowski's integral inequality and the fact that
\[
\norm{|v(s,\cdot)|^{\gamma-1}v(s,\cdot)}_1\leq s^{-\frac{\alpha d(\gamma-1)}{\beta p}}\norm{v}_{E_T}^{\gamma-1}\norm{v(s,\cdot)}_1,\quad s>0,
\]
produce
\[
\left\|\int_0^t Y(t-s,\cdot)\star|v(s,\cdot)|^{\gamma-1}v(s,\cdot)ds\right\|_1\leq \norm{v}_{E_T}^{\gamma}\int_0^t\norm{Y(t-s,\cdot)}_1 s^{-\frac{\alpha d(\gamma-1)}{\beta p}}ds.
\]
From Theorem \ref{CotaNormapY} it follows that
\[
\left\|\int_0^t Y(t-s,\cdot)\star|v(s,\cdot)|^{\gamma-1}v(s,\cdot)ds\right\|_1\lesssim \norm{v}_{E_T}^{\gamma}\int_0^t (t-s)^{\alpha-1} s^{-\frac{\alpha d(\gamma-1)}{\beta p}}ds.
\]
Condition $\frac{d(\gamma-1)}{\beta}<p<\infty$ yields
\[
\left\|\int_0^t Y(t-s,\cdot)\star|v(s,\cdot)|^{\gamma-1}v(s,\cdot)ds\right\|_1\lesssim \norm{v}_{E_T}^{\gamma}t^{\alpha-\frac{\alpha d(\gamma-1)}{\beta p}}
\]
with
\[
0<\alpha-\frac{\alpha d(\gamma-1)}{\beta p}<1.
\]
Similarly,
\[
\left\|\int_0^t Y(t-s,\cdot)\star|v(s,\cdot)|^{\gamma-1}v(s,\cdot)ds\right\|_p\lesssim \norm{v}_{E_T}^{\gamma}t^{\alpha-\frac{\alpha d(\gamma-1)}{\beta p}}.
\]
Therefore, 
\begin{align*}
&\left\|\int_0^t Y(t-s,\cdot)\star|v(s,\cdot)|^{\gamma-1}v(s,\cdot)ds\right\|_1+\left\|\int_0^t Y(t-s,\cdot)\star|v(s,\cdot)|^{\gamma-1}v(s,\cdot)ds\right\|_p\\
&\lesssim \left\|v\right\|_{E_T}^{\gamma}t^{\alpha-\frac{\alpha d(\gamma-1)}{\beta p}}.
\end{align*}
This shows that
\[
\left\|\int_0^t Y(t-s,\cdot)\star|v(s,\cdot)|^{\gamma-1}v(s,\cdot)ds\right\|_{L_p(\RR^d) \cap L_1(\RR^d)}\rightarrow 0
\]
whenever $t\rightarrow 0$, and that
\[
\displaystyle\int_0^tY(t-s,\cdot)\star|v(s,\cdot)|^{\gamma-1}v(s,\cdot)ds \in L_p(\RR^d)\cap L_1(\RR^d),\quad t\in [0,T]. 
\]
Again, condition $\max(1,\kappa) <p<\infty$ and Young's convolution inequality imply that for any $t\in (0,T]$, $Y(t,\cdot)\in L_q(\RR^d)$ with $\frac{1}{p}+\frac{1}{q}=1$ and $1<q<\kappa_2$ (see Theorem \ref{CotaNormapY} and note that $\kappa_2\geq\kappa_1$). Hence, $Y(t-s,\cdot)\star |v(s,\cdot)|^{\gamma-1}v(s,\cdot)\in L_{\infty}(\RR^d)$, $0<s<t$.

\medbreak

Further, condition $\max(1,\kappa)<p<\infty$ also implies that $0<\frac{d}{\beta p}<1$. Consequently,
\begin{align*}
&\left\|\int_0^t Y(t-s,\cdot)\star|v(s,\cdot)|^{\gamma-1}v(s,\cdot)ds\right\|_\infty\\
&\leq \int_0^t\norm{Y(t-s,\cdot)}_{\frac{p}{p-1}} \norm{|v(s,\cdot)|^{\gamma-1}v(s,\cdot)}_p\; ds\\
&\lesssim\norm{v}_{E_T}^{\gamma}\int_0^t(t-s)^{-\frac{\alpha d}{\beta p}+\alpha-1}s^{-\frac{\alpha d(\gamma-1)}{\beta p}}ds\\
&\lesssim\norm{v}_{E_T}^{\gamma}\left[\int_0^{\frac{t}{2}}(t-s)^{-\frac{\alpha d}{\beta p}+\alpha-1}s^{-\frac{\alpha d(\gamma-1)}{\beta p}}ds+\int_{\frac{t}{2}}^t(t-s)^{-\frac{\alpha d}{\beta p}+\alpha-1}s^{-\frac{\alpha d(\gamma-1)}{\beta p}}ds\right]\\
&\lesssim\norm{v}_{E_T}^{\gamma}\left[\left(\frac{t}{2}\right)^{-\frac{\alpha d}{\beta p}+\alpha-1}\int_0^{\frac{t}{2}}s^{-\frac{\alpha d(\gamma-1)}{\beta p}}ds+\left(\frac{t}{2}\right)^{-\frac{\alpha d(\gamma-1)}{\beta p}}\int_{\frac{t}{2}}^t(t-s)^{-\frac{\alpha d}{\beta p}+\alpha-1}ds\right]\\
&\lesssim\norm{v}_{E_T}^{\gamma}\left(\frac{t}{2}\right)^{-\frac{\alpha d}{\beta p}+\alpha-\frac{\alpha d(\gamma-1)}{\beta p}}\\
&<\infty.
\end{align*}
The continuity of
\[
t\mapsto \displaystyle\int_0^t Y(t-s,\cdot)\star|v(s,\cdot)|^{\gamma-1}v(s,\cdot)ds
\]
in $(0,T]$ with respect to the norm topology on $L_p(\RR^d) \cap L_1(\RR^d)$, it follows from conditions on $p$, continuity of $v$ and the property
\begin{equation*}
|\,|a|^c a-|b|^c b\,|\lesssim |\,a-b\,|(|a|^c+|b|^c)\lesssim |\,a-b\,|(\,|a|+|b|\,)^c,\quad a, b\in\RR, c>0.
\end{equation*}
Indeed, suppose $0<t_0<t$ without loss of generality. We have that
\begin{align*}
&\left\|\displaystyle\int_0^t Y(t-s,\cdot)\star|v(s,\cdot)|^{\gamma-1}v(s,\cdot)ds-\displaystyle\int_0^{t_0} Y(t_0-s,\cdot)\star|v(s,\cdot)|^{\gamma-1}v(s,\cdot)ds\right\|_1\\
&\leq\displaystyle\int_0^{t_0} \norm{Y(s,\cdot)}_1\norm{|v(t-s,\cdot)|^{\gamma-1}v(t-s,\cdot)-|v(t_0-s,\cdot)|^{\gamma-1}v(t_0-s,\cdot)}_1 ds\\
&~~~~+\displaystyle\int_{t_0}^t \norm{Y(s,\cdot)}_1\norm{|v(t-s,\cdot)|^{\gamma-1}v(t-s,\cdot)}_1ds\\
&\lesssim \norm{v}_{E_T}^{\gamma-1}\varepsilon\displaystyle\int_0^{t_0} s^{\alpha-1}\left((t-s)^{-\frac{\alpha d(\gamma-1)}{\beta p}}+(t_0-s)^{-\frac{\alpha d(\gamma-1)}{\beta p}}\right) ds\\
&~~~~+\norm{v}_{E_T}^{\gamma}\displaystyle\int_{t_0}^t s^{\alpha-1}(t-s)^{-\frac{\alpha d(\gamma-1)}{\beta p}} ds\\
&\lesssim \norm{v}_{E_T}^{\gamma-1}\varepsilon\displaystyle\int_0^{t_0} s^{\alpha-1}(t_0-s)^{-\frac{\alpha d(\gamma-1)}{\beta p}}ds+\norm{v}_{E_T}^{\gamma}\displaystyle\int_{t_0}^t s^{\alpha-1}(t-s)^{-\frac{\alpha d(\gamma-1)}{\beta p}} ds\\
&\lesssim \norm{v}_{E_T}^{\gamma-1}\varepsilon t_0^{\alpha-\frac{\alpha d(\gamma-1)}{\beta p}}+\norm{v}_{E_T}^{\gamma}(t-t_0)^{\alpha-\frac{\alpha d(\gamma-1)}{\beta p}}.
\end{align*}
Similarly, we obtain that
\[
\left\|\displaystyle\int_0^t Y(t-s,\cdot)\star|v(s,\cdot)|^{\gamma-1}v(s,\cdot)ds-\displaystyle\int_0^{t_0} Y(t_0-s,\cdot)\star|v(s,\cdot)|^{\gamma-1}v(s,\cdot)ds\right\|_p\rightarrow 0
\]
whenever $t\rightarrow t_0$.

Now, for the continuity of
\[
t\mapsto \displaystyle\int_0^t Y(t-s,\cdot)\star|v(s,\cdot)|^{\gamma-1}v(s,\cdot)ds
\]
in $(0,T]$ with respect to the norm topology on $L_{\infty}(\RR^d)$, we find that 
\begin{align*}
&\left\|\displaystyle\int_0^t Y(t-s,\cdot)\star|v(s,\cdot)|^{\gamma-1}v(s,\cdot)ds-\displaystyle\int_0^{t_0} Y(t_0-s,\cdot)\star|v(s,\cdot)|^{\gamma-1}v(s,\cdot)ds\right\|_\infty\\
&\leq\displaystyle\int_0^{t_0} \norm{Y(s,\cdot)}_q\norm{|v(t-s,\cdot)|^{\gamma-1}v(t-s,\cdot)-|v(t_0-s,\cdot)|^{\gamma-1}v(t_0-s,\cdot)}_p \;ds\\
&~~~~+\displaystyle\int_{t_0}^t \norm{Y(s,\cdot)}_q\norm{|v(t-s,\cdot)|^{\gamma-1}v(t-s,\cdot)}_p\;ds\\
&\lesssim \norm{v}_{E_T}^{\gamma-1}\varepsilon\displaystyle\int_0^{t_0} s^{-\frac{\alpha d}{\beta p}+\alpha-1}\left((t-s)^{-\frac{\alpha d(\gamma-1)}{\beta p}}+(t_0-s)^{-\frac{\alpha d(\gamma-1)}{\beta p}}\right) ds\\
&~~~~ +\norm{v}_{E_T}^{\gamma}\displaystyle\int_{t_0}^t s^{-\frac{\alpha d}{\beta p}+\alpha-1}(t-s)^{-\frac{\alpha d(\gamma-1)}{\beta p}} ds\\
&\lesssim \norm{v}_{E_T}^{\gamma-1}\varepsilon\displaystyle\int_0^{t_0} s^{-\frac{\alpha d}{\beta p}+\alpha-1}(t_0-s)^{-\frac{\alpha d(\gamma-1)}{\beta p}}ds+\norm{v}_{E_T}^{\gamma}\displaystyle\int_{t_0}^t s^{-\frac{\alpha d}{\beta p}+\alpha-1}(t-s)^{-\frac{\alpha d(\gamma-1)}{\beta p}} ds\\
&\lesssim \norm{v}_{E_T}^{\gamma-1}\varepsilon t_0^{-\frac{\alpha d}{\beta p}+\alpha-\frac{\alpha d(\gamma-1)}{\beta p}}+\norm{v}_{E_T}^{\gamma}t_0^{-\frac{\alpha d}{\beta p}}(t-t_0)^{\alpha-\frac{\alpha d(\gamma-1)}{\beta p}}.
\end{align*}
Up to this point, we have proved that the operator $\mathcal{M}$ given by \eqref{operator} is well defined.

Now, let $v,w\in E_{T}$. Previous arguments show that
\begin{equation*}
\begin{split}
&\left\|\mathcal{M}(v)(t,\cdot)-\mathcal{M}(w)(t,\cdot)\right\|_1\\
&\leq|\lambda|\displaystyle\int_0^t \norm{Y(t-s,\cdot)}_1 \left\||v(s,\cdot)|^{\gamma-1}v(s,\cdot)-|w(s,\cdot)|^{\gamma-1}w(s,\cdot)\right\|_1 ds\\
&\lesssim \left(\norm{v}_{E_T}+\norm{w}_{E_T}\right)^{\gamma-1}\norm{v-w}_{E_T}\displaystyle\int_0^{t} (t-s)^{\alpha-1}s^{-\frac{\alpha d(\gamma-1)}{\beta p}}ds\\
&\lesssim \left(\norm{v}_{E_T}+\norm{w}_{E_T}\right)^{\gamma-1}\norm{v-w}_{E_T}t^{\alpha-\frac{\alpha d(\gamma-1)}{\beta p}}\\
&\lesssim \left(\norm{v}_{E_T}+\norm{w}_{E_T}\right)^{\gamma-1}\norm{v-w}_{E_T}T^{\alpha-\frac{\alpha d(\gamma-1)}{\beta p}}.
\end{split}
\end{equation*}
In the same way we estimate
\[
\left\|\mathcal{M}(v)(t,\cdot)-\mathcal{M}(w)(t,\cdot)\right\|_p\lesssim \left(\norm{v}_{E_T}+\norm{w}_{E_T}\right)^{\gamma-1}\norm{v-w}_{E_T}T^{\alpha-\frac{\alpha d(\gamma-1)}{\beta p}}.
\]
Similarly,
\[
\left\|\mathcal{M}(v)(t,\cdot)-\mathcal{M}(w)(t,\cdot)\right\|_\infty\lesssim \left(\norm{v}_{E_T}+\norm{w}_{E_T}\right)^{\gamma-1}\norm{v-w}_{E_T}t^{\alpha-\frac{\alpha d}{\beta p}-\frac{\alpha d(\gamma-1)}{\beta p}}
\]
and we have
\begin{equation*}
\begin{split}
t^{\frac{\alpha d}{\beta p}}\left\|\mathcal{M}(v)(t,\cdot)-\mathcal{M}(w)(t,\cdot)\right\|_\infty&\lesssim \left(\norm{v}_{E_T}+\norm{w}_{E_T}\right)^{\gamma-1}\norm{v-w}_{E_T}t^{\alpha-\frac{\alpha d(\gamma-1)}{\beta p}}\\
&\lesssim \left(\norm{v}_{E_T}+\norm{w}_{E_T}\right)^{\gamma-1}\norm{v-w}_{E_T}T^{\alpha-\frac{\alpha d(\gamma-1)}{\beta p}}.
\end{split}
\end{equation*}
This shows that
\begin{equation}
\label{cotapartenolineal}
\left\|\mathcal{M}(v)-\mathcal{M}(w)\right\|_{E_T}\leq C_2T^{\alpha-\frac{\alpha d(\gamma-1)}{\beta p}}\norm{v-w}_{E_T}\left(\norm{v}_{E_T}+\norm{w}_{E_T}\right)^{\gamma-1}.
\end{equation}
Besides, we derive from \eqref{cotasiniciales} that
\begin{equation}
\label{cotapartelineal}
\norm{Z\star u_0}_{E_T}\leq C_1(\norm{u_0}_1+\norm{u_0}_p).
\end{equation}
Let $T^*\in(0,T)$ and $R=2C_1\norm{u_0}_{L_p(\RR^d) \cap L_1(\RR^d)}$. We consider the closed ball
\[
B_{T^*,R}:=\{w\in E_{T^*}:\norm{w}_{E_{T^*}}\leq R\}.
\]
Our aim now is to obtain a suitable $T^*$ such that $\mathcal{M}$ is a contraction as an operator $B_{T^*,R}\rightarrow B_{T^*,R}$, thereby existence and uniqueness of a fixed point of this operator follow from Banach fixed-point theorem.

First, we find a condition on $T^*$ such that $\mathcal{M}(B_{T^*,R})\subset B_{T^*,R}$. Let $w\in B_{T^*,R}$. Using \eqref{cotapartenolineal} with $v\equiv 0$ and \eqref{cotapartelineal}, we obtain that
\begin{align*}
\norm{\mathcal{M}(w)}_{E_{T^*}}&\leq C_1\norm{u_0}_{L_p(\RR^d) \cap L_1(\RR^d)}+C_2(T^*)^{\alpha-\frac{\alpha d(\gamma-1)}{\beta p}}\norm{w}_{E_{T^*}}\left(\norm{w}_{E_{T^*}}\right)^{\gamma-1}\\
&=\dfrac{R}{2}+C_2(T^*)^{\alpha-\frac{\alpha d(\gamma-1)}{\beta p}}\norm{w}_{E_{T^*}}^{\gamma}\\
&\leq \dfrac{R}{2}+C_2(T^*)^{\alpha-\frac{\alpha d(\gamma-1)}{\beta p}}R^{\gamma}.
\end{align*}
Therefore, we need to set the condition 
\begin{equation}
\label{Condicion1T^*}
C_2(T^*)^{\alpha-\frac{\alpha d(\gamma-1)}{\beta p}}R^{\gamma}\leq \dfrac{R}{2}. 
\end{equation}
On the other hand, let $v,w\in B_{T^*,R}$. Using \eqref{cotapartenolineal} we get
\begin{align*}
\norm{\mathcal{M}(v)-\mathcal{M}(w)}_{E_{T^*}}&\leq C_2(T^*)^{\alpha-\frac{\alpha d(\gamma-1)}{\beta p}}\norm{v-w}_{E_{T^*}}\left(\norm{v}_{E_{T^*}}+\norm{w}_{E_{T^*}}\right)^{\gamma-1}\\
&\leq C_2(T^*)^{\alpha-\frac{\alpha d(\gamma-1)}{\beta p}}\norm{v-w}_{E_{T^*}}(2R)^{\gamma-1}.
\end{align*}
Consequently, we also need 
\begin{equation}
\label{Condicion2T^*}
C_2(T^*)^{\alpha-\frac{\alpha d(\gamma-1)}{\beta p}}(2R)^{\gamma-1} < 1. 
\end{equation}
Thus, for sufficiently small $T^*$, the requirements \eqref{Condicion1T^*} and \eqref{Condicion2T^*} are satisfied.
\end{proof}

\section{Global solution for small initial condition}

\label{sec:3}
By a global solution of the Cauchy problem \eqref{general} we understand the solution $u$ of the integral equation \eqref{Inteq:1} belonging to the Banach space
\[
E:=C([0,\infty); L_p(\RR^d) \cap L_1(\RR^d)) \cap C((0,\infty); L_{\infty}(\RR^d)),
\]
with the norm
	\[
	\left\|v\right\|_{E}:=\sup_{t\geq 0}\left(\seq{t}^{\frac{\alpha d}{\beta}\left(\frac{1}{p'}-\frac{1}{p}\right)}\left\|v(t,\cdot)\right\|_p + \left\|v(t,\cdot)\right\|_1\right) + \sup_{t>0}\{t\}^{\frac{\alpha d}{\beta p}}\seq{t}^{\frac{\alpha d}{\beta p'}}\left\|v(t,\cdot)\right\|_{\infty},
	\]
where $1\leq p'<p$, $\seq{t}:=\sqrt{1+t^2}$ and $\{t\}:=\dfrac{t}{\sqrt{1+t^2}}$.

\medbreak

As in the previous section, we define on $E$ the operator  
\[
\mathcal{M}(v)(t,x):=\int_{\RR^d} Z(t,x-y)u_0(y)dy + \lambda\int_0^t\int_{\RR^d} Y(t-s,x-y)|v(s,y)|^{\gamma-1}v(s,y)dyds,
\]
$u_0\in L_p(\RR^d) \cap L_1(\RR^d)$ is a given data and $v\in E$.

Let $\lambda\in \RR$ and $\gamma>1$. Let $1=p'<\frac{d}{\beta}(\gamma-1)$ whenever $d<\beta$, or $\frac{d}{\beta}<p'<\frac{d}{\beta}(\gamma-1)$ whenever $d\geq\beta$. Suppose that $\max\left(1,\kappa,\frac{d(\gamma-1)}{\beta}\right)<p<\infty$. These conditions on $p$ and $p'$ guarantee the existence of $1\leq q<\kappa_1$ such that
\[
\frac{1}{l}+\frac{1}{q}=1+\frac{1}{r},
\]
considering $r\in\{1,p,\infty\}$, where $l$, $q$ and $r$ are related to Young's convolution inequality, that is,
\[
\norm{Z(t,\cdot)\star u_0}_r\leq\norm{Z(t,\cdot)}_q\norm{u_0}_l,\quad t\geq 0.
\]
We note that for $r=1$, this is possible only if $l=q=1$. Similar situation we have for $Y$. For $r=1$ we can't get a factor of time estimating $\norm{Z\star u_0}_1$ because $\norm{Z(t,\cdot)}_1=1$ por all $t\geq 0$.

These conditions, together with Young's convolution inequality and interpolation, yield the following bounds.

$\left\|Z(t,\cdot)\star u_0\right\|_1\leq\left\|Z(t,\cdot)\right\|_1\left\|u_0\right\|_1=\left\|u_0\right\|_1$ for all $t\geq 0$.

$\left\|Z(t,\cdot)\star u_0\right\|_p\leq\left\|Z(t,\cdot)\right\|_1\left\|u_0\right\|_p=\left\|u_0\right\|_p$ for all $0\leq t\leq 1$.

$\left\|Z(t,\cdot)\star u_0\right\|_p\leq\left\|Z(t,\cdot)\right\|_q\left\|u_0\right\|_{p'}\lesssim t^{-\frac{\alpha d}{\beta}\left(\frac{1}{p'}-\frac{1}{p}\right)}\max (\norm{u_0}_1,\norm{u_0}_p)$ for all $t> 1$.

$\left\|Z(t,\cdot)\star u_0\right\|_{\infty}\leq\left\|Z(t,\cdot)\right\|_q\left\|u_0\right\|_{p}\lesssim t^{-\frac{\alpha d}{\beta p}}\left\|u_0\right\|_p$ for all $0<t\leq 1$.

$\left\|Z(t,\cdot)\star u_0\right\|_{\infty}\leq\left\|Z(t,\cdot)\right\|_q\left\|u_0\right\|_{p'}\lesssim t^{-\frac{\alpha d}{\beta p'}}\left\|u_0\right\|_{p'}\lesssim t^{-\frac{\alpha d}{\beta p'}}\max (\norm{u_0}_1,\norm{u_0}_p)$ for $t> 1$.

Therefore,
\begin{align*}
\sup_{t\geq 0}\left(\langle t\rangle^{\frac{\alpha d}{\beta}\left(\frac{1}{p'}-\frac{1}{p}\right)}\left\|Z(t,\cdot)\star u_0\right\|_p+\left\|Z(t,\cdot)\star u_0\right\|_1\right)+
\sup_{t>0}\{t\}^{\frac{\alpha d}{\beta p}}\langle t\rangle^{\frac{\alpha d}{\beta p'}}\left\|Z(t,\cdot)\star u_0\right\|_{\infty}\lesssim (\norm{u_0}_1+\norm{u_0}_p).
\end{align*}
For $0\leq t\leq 1$ we have that
\begin{equation*}
\begin{split}
\left\|\displaystyle\int_0^t Y(t-s,\cdot)\star|v(s,\cdot)|^{\gamma-1}v(s,\cdot)ds\right\|_1&\leq \displaystyle\int_0^{t} \norm{Y(t-s,\cdot)}_1\norm{|v(s,\cdot)|^{\gamma-1}v(s,\cdot)}_1 ds\\
&\lesssim \norm{v}_E^\gamma\int_0^t(t-s)^{\alpha-1}s^{-\frac{\alpha d}{\beta p}(\gamma-1)}\langle s\rangle^{-\frac{\alpha d}{\beta}(\gamma-1)\left(\frac{1}{p'}-\frac{1}{p}\right)}ds\\
&\lesssim \norm{v}_E^\gamma\int_0^t(t-s)^{\alpha-1}s^{-\frac{\alpha d}{\beta p}(\gamma-1)}ds\\
&\lesssim \norm{v}_E^\gamma t^{\alpha-\frac{\alpha d}{\beta p}(\gamma-1)}\lesssim \norm{v}_E^\gamma.
\end{split}
\end{equation*}
For all $t>1$, we use the fact that $\alpha-\frac{\alpha d(\gamma-1)}{\beta p}>0$ and that $\alpha-\frac{\alpha d(\gamma-1)}{\beta p'}<0$. Choosing $0<c<\dfrac{t}{2}$ we get
\begin{align*}
&\left\|\displaystyle\int_0^t Y(t-s,\cdot)\star|v(s,\cdot)|^{\gamma-1}v(s,\cdot)ds\right\|_1\\
&\leq \displaystyle\int_0^{t} \norm{Y(t-s,\cdot)}_{1}\norm{|v(s,\cdot)|^{\gamma-1}v(s,\cdot)}_{1} ds\\
&\lesssim \norm{v}_E^{\gamma}\int_0^t(t-s)^{\alpha-1}s^{-\frac{\alpha d}{\beta p}(\gamma-1)}\langle s\rangle^{-\frac{\alpha d}{\beta}(\gamma-1)\left(\frac{1}{p'}-\frac{1}{p}\right)} ds\\
&\lesssim\norm{v}_E^{\gamma}\left[\displaystyle\int_0^{\frac{t}{2}} (t-s)^{\alpha-1}s^{-\frac{\alpha d}{\beta p}(\gamma-1)}\langle s\rangle^{-\frac{\alpha d}{\beta}(\gamma-1)\left(\frac{1}{p'}-\frac{1}{p}\right)} ds\right.\\
&~~~~\left.+\int_{\frac{t}{2}}^t (t-s)^{\alpha-1}s^{-\frac{\alpha d}{\beta p}(\gamma-1)}\langle s\rangle^{-\frac{\alpha d}{\beta}(\gamma-1)\left(\frac{1}{p'}-\frac{1}{p}\right)} ds\right]\\
&\lesssim \norm{v}_E^{\gamma}\left[\displaystyle\int_0^{c}s^{\alpha-1-\frac{\alpha d}{\beta p}(\gamma-1)}\langle s\rangle^{-\frac{\alpha d}{\beta}(\gamma-1)\left(\frac{1}{p'}-\frac{1}{p}\right)} ds\right.\\
&~~~~+\displaystyle\int_c^{\frac{t}{2}}s^{\alpha-1-\frac{\alpha d}{\beta p}(\gamma-1)}\langle s\rangle^{-\frac{\alpha d}{\beta}(\gamma-1)\left(\frac{1}{p'}-\frac{1}{p}\right)} ds\\
&~~~~\left.+\int_{\frac{t}{2}}^t (t-s)^{\alpha-1}s^{-\frac{\alpha d}{\beta p}(\gamma-1)}\langle s\rangle^{-\frac{\alpha d}{\beta}(\gamma-1)\left(\frac{1}{p'}-\frac{1}{p}\right)} ds\right]\\
&\lesssim \norm{v}_E^{\gamma}\left[\displaystyle\int_0^{c}s^{\alpha-1-\frac{\alpha d}{\beta p}(\gamma-1)} ds+\displaystyle\int_c^{\frac{t}{2}}s^{\alpha-1-\frac{\alpha d}{\beta p'}(\gamma-1)} ds\right.\\
&~~~~\left. +\left\langle\frac{t}{2}\right\rangle^{-\frac{\alpha d}{\beta}(\gamma-1)\left(\frac{1}{p'}-\frac{1}{p}\right)}\displaystyle\int_{\frac{t}{2}}^t (t-s)^{\alpha-1}s^{-\frac{\alpha d}{\beta p}(\gamma-1)}ds\right]\\
&\lesssim \norm{v}_E^{\gamma}\left[\frac{c^{\alpha-\frac{\alpha d}{\beta p}(\gamma-1)}}{\alpha-\frac{\alpha d(\gamma-1)}{\beta p}}-\frac{c^{\alpha-\frac{\alpha d}{\beta p'}(\gamma-1)}}{\alpha-\frac{\alpha d(\gamma-1)}{\beta p'}}+\left\langle\frac{t}{2}\right\rangle^{-\frac{\alpha d}{\beta}(\gamma-1)\left(\frac{1}{p'}-\frac{1}{p}\right)}\left\langle\frac{t}{2}\right\rangle^{\alpha-\frac{\alpha d}{\beta p}(\gamma-1)}\right]\\
&\lesssim \norm{v}_E^{\gamma}\left[C+\left\langle\frac{t}{2}\right\rangle^{\alpha-\frac{\alpha d}{\beta p'}(\gamma-1)}\right]\lesssim \norm{v}_E^\gamma.
\end{align*}
For $0\leq t\leq 1$ we get
\begin{equation*}
\begin{split}
&\left\|\displaystyle\int_0^t Y(t-s,\cdot)\star|v(s,\cdot)|^{\gamma-1}v(s,\cdot)ds\right\|_p\\
&\leq \displaystyle\int_0^{t} \norm{Y(t-s,\cdot)}_1\norm{|v(s,\cdot)|^{\gamma-1}v(s,\cdot)}_p ds\\
&\lesssim \norm{v}_E^\gamma\int_0^t(t-s)^{\alpha-1}s^{-\frac{\alpha d}{\beta p}(\gamma-1)}\langle s\rangle^{-\frac{\alpha d}{\beta}(\gamma-1)\left(\frac{1}{p'}-\frac{1}{p}\right)} \langle s\rangle^{-\frac{\alpha d}{\beta}\left(\frac{1}{p'}-\frac{1}{p}\right)} ds\\
&\lesssim \norm{v}_E^\gamma\int_0^t(t-s)^{\alpha-1}s^{-\frac{\alpha d}{\beta p}(\gamma-1)}ds\\
&\lesssim\norm{v}_E^\gamma.
\end{split}
\end{equation*}
For $t > 1$, the fact that $\norm{v(s,\cdot)}_{p'}\leq \max(\norm{v(s,\cdot)}_{1},\norm{v(s,\cdot)}_{p})$, $s>0$, together with the condition $\frac{d}{\beta p'}<1$, yields 
\begin{equation*}
\begin{split}
&\left\|\displaystyle\int_0^t Y(t-s,\cdot)\star|v(s,\cdot)|^{\gamma-1}v(s,\cdot)ds\right\|_p\\
&\leq \displaystyle\int_0^{t} \norm{Y(t-s,\cdot)}_{q}\norm{|v(s,\cdot)|^{\gamma-1}v(s,\cdot)}_{p'} ds\\
&\lesssim \norm{v}_E^{\gamma-1}\int_0^t(t-s)^{-\frac{\alpha d}{\beta}\left(\frac{1}{p'}-\frac{1}{p}\right)+\alpha-1}s^{-\frac{\alpha d}{\beta p}(\gamma-1)}\langle s\rangle^{-\frac{\alpha d}{\beta}(\gamma-1)\left(\frac{1}{p'}-\frac{1}{p}\right)} \norm{v(s,\cdot)}_{p'} ds\\
&\lesssim \norm{v}_E^{\gamma}\int_0^t(t-s)^{-\frac{\alpha d}{\beta}\left(\frac{1}{p'}-\frac{1}{p}\right)+\alpha-1}s^{-\frac{\alpha d}{\beta p}(\gamma-1)}\langle s\rangle^{-\frac{\alpha d}{\beta}(\gamma-1)\left(\frac{1}{p'}-\frac{1}{p}\right)}ds\\
&\lesssim \norm{v}_E^{\gamma}t^{\frac{\alpha d}{\beta p}}\int_0^t(t-s)^{-\frac{\alpha d}{\beta p'}+\alpha-1}s^{-\frac{\alpha d}{\beta p}(\gamma-1)}\langle s\rangle^{-\frac{\alpha d}{\beta}(\gamma-1)\left(\frac{1}{p'}-\frac{1}{p}\right)}ds\\
&\lesssim \norm{v}_E^{\gamma}t^{\frac{\alpha d}{\beta p}}\left[\int_0^{\frac{t}{2}}(t-s)^{-\frac{\alpha d}{\beta p'}+\alpha-1}s^{-\frac{\alpha d}{\beta p}(\gamma-1)}\langle s\rangle^{-\frac{\alpha d}{\beta}(\gamma-1)\left(\frac{1}{p'}-\frac{1}{p}\right)}ds\right.\\
&\left.+\int_{\frac{t}{2}}^t(t-s)^{-\frac{\alpha d}{\beta p'}+\alpha-1}s^{-\frac{\alpha d}{\beta p}(\gamma-1)}\langle s\rangle^{-\frac{\alpha d}{\beta}(\gamma-1)\left(\frac{1}{p'}-\frac{1}{p}\right)}ds\right]\\
&\leq \norm{v}_E^{\gamma}t^{\frac{\alpha d}{\beta p}}\left[\left(\dfrac{t}{2}\right)^{-\frac{\alpha d}{\beta p'}}\int_0^{\frac{t}{2}}(t-s)^{\alpha-1}s^{-\frac{\alpha d}{\beta p}(\gamma-1)}\langle s\rangle^{-\frac{\alpha d}{\beta}(\gamma-1)\left(\frac{1}{p'}-\frac{1}{p}\right)}ds\right.\\
&\left.+\left\langle\frac{t}{2}\right\rangle^{-\frac{\alpha d}{\beta}(\gamma-1)\left(\frac{1}{p'}-\frac{1}{p}\right)}\left(\dfrac{t}{2}\right)^{-\frac{\alpha d}{\beta p}(\gamma-1)}\int_{\frac{t}{2}}^t(t-s)^{-\frac{\alpha d}{\beta p'}+\alpha-1}ds\right]\\
&\lesssim \norm{v}_E^{\gamma}t^{\frac{\alpha d}{\beta p}}\left[C\left(\dfrac{t}{2}\right)^{-\frac{\alpha d}{\beta p'}}+\left\langle\frac{t}{2}\right\rangle^{-\frac{\alpha d}{\beta}(\gamma-1)\left(\frac{1}{p'}-\frac{1}{p}\right)}\left(\dfrac{t}{2}\right)^{-\frac{\alpha d}{\beta p}(\gamma-1)-\frac{\alpha d}{\beta p'}+\alpha}\right]\\
&\lesssim \norm{v}_E^{\gamma}t^{\frac{\alpha d}{\beta p}-\frac{\alpha d}{\beta p'}}\left[C+\left\langle\frac{t}{2}\right\rangle^{-\frac{\alpha d}{\beta}(\gamma-1)\left(\frac{1}{p'}-\frac{1}{p}\right)}\left(\dfrac{t}{2}\right)^{-\frac{\alpha d}{\beta p}(\gamma-1)+\alpha}\right]\\
&\lesssim \norm{v}_E^\gamma t^{-\frac{\alpha d}{\beta}\left(\frac{1}{p'}-\frac{1}{p}\right)}.
\end{split}
\end{equation*}
Now, for all $0<t\leq 1$, as in the proof of Theorem \ref{localsolution}, we obtain that
\begin{equation*}
\begin{split}
&\left\|\displaystyle\int_0^t Y(t-s,\cdot)\star|v(s,\cdot)|^{\gamma-1}v(s,\cdot)ds\right\|_{\infty}\\
&\leq \displaystyle\int_0^{t} \norm{Y(t-s,\cdot)}_{\frac{p}{p-1}}\norm{|v(s,\cdot)|^{\gamma-1}v(s,\cdot)}_p ds\\
&\lesssim\norm{v}_E^\gamma \int_0^t (t-s)^{-\frac{\alpha d}{\beta p}+\alpha-1}s^{-\frac{\alpha d}{\beta p}(\gamma-1)}\langle s\rangle^{-\frac{\alpha d}{\beta}(\gamma-1)\left(\frac{1}{p'}-\frac{1}{p}\right)}ds\\
&\lesssim\norm{v}_E^\gamma \int_0^t (t-s)^{-\frac{\alpha d}{\beta p}+\alpha-1}s^{-\frac{\alpha d}{\beta p}(\gamma-1)}ds\\
&\lesssim\norm{v}_E^\gamma t^{-\frac{\alpha d}{\beta p}}.
\end{split}
\end{equation*}
For all $t>1$ we find
\begin{align*} 
&\left\|\displaystyle\int_0^t Y(t-s,\cdot)\star|v(s,\cdot)|^{\gamma-1}v(s,\cdot)ds\right\|_\infty\\
&\leq\displaystyle\int_0^{t} \norm{Y(t-s,\cdot)}_{\frac{p'}{p'-1}}\norm{|v(s,\cdot)|^{\gamma-1}v(s,\cdot)}_{p'} ds\\
&\lesssim \norm{v}_E^\gamma\displaystyle\int_0^{t}(t-s)^{-\frac{\alpha d}{\beta p'}+\alpha-1}s^{-\frac{\alpha d}{\beta p}(\gamma-1)}\langle s\rangle^{-\frac{\alpha d}{\beta}(\gamma-1)\left(\frac{1}{p'}-\frac{1}{p}\right)}ds\\
&\lesssim \norm{v}_E^\gamma t^{-\frac{\alpha d}{\beta p'}}.
\end{align*}
Using these bounds and from a straightforward inspection of the proof of Theorem \ref{localsolution}, it follows that the operator $\mathcal{M}$ is well defined on $E$. From the beginning of this section, we also get the estimate
\[
\norm{Z\star u_0}_{E}\leq C_1(\norm{u_0}_1+\norm{u_0}_p).
\]
In the same way, as in the proof of Theorem \ref{localsolution}, the following estimate
\[
\left\|\mathcal{M}(v)-\mathcal{M}(w)\right\|_{E}\leq C_2\norm{v-w}_{E}\left(\norm{v}_{E}+\norm{w}_{E}\right)^{\gamma-1}
\]
holds, considering as before the cases $0\leq t\leq 1$ and $t>1$, respectively.

Finally, if $u_0\in L_p(\RR^d) \cap L_1(\RR^d)$ is sufficiently small, then the operator $\mathcal{M}$ also satisfies the proof of Theorem \ref{localsolution} with the corresponding closed ball on $E$, that is, the conditions \eqref{Condicion1T^*} and \eqref{Condicion2T^*} are satisfied without the restriction on the time.

Consequently, we have proved the following result.

\begin{theorem}
\label{globalsolutionsmallu0}
Let $\alpha\in (0,1)$ and $\beta\in (0,2)$. Assume the hypothesis $(\mathcal{H}_1)$ holds. Let $\lambda\in \RR$ and $\gamma>1$. Let $1=p'<\frac{d}{\beta}(\gamma-1)$ whenever $d<\beta$, or $\frac{d}{\beta}<p'<\frac{d}{\beta}(\gamma-1)$ whenever $d\geq\beta$. Suppose that $\max\left(1,\kappa,\frac{d(\gamma-1)}{\beta}\right)<p<\infty$. If $u_0\in L_p(\RR^d) \cap L_1(\RR^d)$ is sufficiently small, then the operator $\mathcal{M}$ has a unique fixed point $u$ in $E$ and the optimal time decay estimate
\[
\norm{u(t)}_1+t^{\frac{\alpha d}{\beta}\left(\frac{1}{p'}-\frac{1}{p}\right)}\norm{u(t)}_p+t^{\frac{\alpha d}{\beta p'}}\norm{u(t)}_\infty\lesssim(\norm{u_0}_1+\norm{u_0}_p) 
\]  
is true for all $t\geq 1$.
\end{theorem}

\section{Global non-negative solution}
\label{sec:4}
In this section we consider the operator $(-\Psi_{\beta}(-i\nabla), C_0^{\infty}(\RR^d))$ on the Hilbert
space $X=(L_2(\RR^d), \norm{\cdot}_2)$. First, we have that $(-\Psi_{\beta}(-i\nabla), C_0^{\infty}(\RR^d))$ satisfies the positive maximum principle since its symbol $\psi$ is a continuous and negative definite function on $\RR^d$ (\cite[Theorem 4.5.6]{Jac01}). It is also symmetric because $\psi$ is real and symmetric. Consequently, $(-\Psi_{\beta}(-i\nabla), C_0^{\infty}(\RR^d))$ is closable (\cite[Theorem 3.6]{Schi01}).

Defining the domain of $-\Psi_{\beta}(-i\nabla)$ as $D(-\Psi_{\beta}(-i\nabla))=\overline{C_0^\infty(\RR^d)}^{\norm{\cdot}_{\Psi_{\beta},L_2}}=:H_2^\beta(\RR^d)$, where the clausure is respect to the graph norm $\norm{\cdot}_{\Psi_{\beta},L_2}^2=\norm{\cdot}_{2}^2 + \norm{\Psi_{\beta}(-i\nabla)(\cdot)}_{2}^2$ (\cite[Theorems 2.7.14 and 3.10.3]{Jac01}) and $H_2^\beta(\RR^d)$ is an \textit{anisotropic function space} (\cite[Section 3.10 and Example 4.1.16]{Jac01}), we get that $(-\Psi_{\beta}(-i\nabla),H_2^\beta(\RR^d))$ generates a symmetric sub-Markovian semigroup on $L_2(\RR^d)$, (\cite[Examples 4.1.13 and 4.3.9]{Jac01}).

By definition \cite[Definition 4.1.6]{Jac01}, a sub-Markovian semigroup is a strongly continuous contraction semigroup. Therefore, $(-\Psi_{\beta}(-i\nabla),H_2^\beta(\RR^d))$ is closed and $H_2^\beta(\RR^d)$ is dense in $L_2(\RR^d)$ (\cite[Corollary 4.1.15]{Jac01}). Moreover, $(-\Psi_{\beta}(-i\nabla),H_2^\beta(\RR^d))$ is a self-adjoint operator (\cite[Section 4.7]{Jac01}).

Since $-\Psi_{\beta}(-i\nabla)=:A$ is closed, linear, densely defined and self-adjoint on the Hilbert space $X$, it follows that $-A$ is a normal operator (\cite[Definition 13.29]{Rud02}) and $\sigma(-A)\subset\mathbb{R}$. Moreover, $\sigma(-A)\subset [0,\infty)$ because $A$ satisfies \cite[Theorem 8.3.2 part (i)]{Dav07}. From Parseval's theorem (or Courrège, see \eqref{ProductoInterno} below) we also have that $-A$ is strictly positive, hence it is 1-1 and satisfies \cite[Theorem 13.11 b)]{Rud02}. Therefore, $-A$ is sectorial (\cite[Section 8.1]{Pru93}). Besides, $X$ belongs to the class $\mathcal{HT}$ (see \cite[definition in page 216, a characterization in page 217 and page 234]{Pru93}).

This shows that the operator $-A$ belongs to $\mathcal{BIP}(X)$ and $\theta_{-A}=0$ (\cite[Definition 8.1 and Section 8.7 c)(i)]{Pru93}), furthermore, it satisfies \cite[Theorem 8.7 part (i)]{Pru93} with $\omega_{-A}=0$.

On the other hand, the Laplace-transform of $g_\alpha(t)=\frac{t^{\alpha-1}}{\Gamma(\alpha)},\;t>0,$ is $\hat{g_\alpha}(s)=s^{-\alpha}$, $\text{Re}(s)>0$ (\cite[Example 2.1]{Pru93}). This yields
\[
\lim_{s\rightarrow\infty}|\hat{g_\alpha}(s)|<\infty.
\]
The kernel $g_\alpha$ is also $1$-regular (\cite[Definition 3.4 and Proposition 3.3]{Pru93}) and $\theta_a$-sectorial with $\theta_a=\alpha\frac{\pi}{2}$ (\cite[Definition 3.2]{Pru93}). Therefore, $\theta_a+\theta_{-A}<\pi$.

This shows that $g_\alpha$ satisfies \cite[Theorem 8.7 parts (ii), (iv) and (v)]{Pru93}, with $\omega_a=0$. 

Now, we consider the Volterra equation
\begin{equation}
\label{Volterrageneral}
u(t)=f(t)+\int_0^t g_\alpha(t-s)Au(s)ds,\quad t\in [0,T].
\end{equation}  
The family of bounded linear operators $\{S(t)\}_{t\geq 0}$ on $L_2(\RR^d)$, given by
\begin{equation}
\label{Resolvente}
S(t)v:=Z(t,\cdot)\star v,
\end{equation}
is a \textbf{resolvent} for \eqref{Volterrageneral}. That is, $S$ satisfies the following conditions (\cite[Definition 1.3]{Pru93}).
\begin{align*}
(S_1)&: S(0)=I \text{ and } S(t)\text{ is strongly continuous on }[0,\infty),\\
(S_2)&: S(t)v\in D(A) \text{ and } AS(t)v =S(t)Av, \text{ for all }v\in D(A)\text{ and }t \geq 0,\\
(S_3)&: S(t)v=v+\int_0^tg_\alpha(t-s)AS(s)vds, \text{ for all }v\in D(A)\text{ and }t \geq 0.
\end{align*}
Indeed, it is easy to see that $S(t)$ is a bounded linear operator and $(S_1)$ is satisfied from Theorem \ref{C1menosZ}. For $(S_2)$, let $t\ge 0$ and $v\in H_2^\beta(\RR^d)$. We have that
\[
\norm{S(t)v}_2\leq\norm{v}_2<\infty
\]
and
\begin{align*}
\int_{\RR^d}(1+\norm{\xi}^2)^\beta|\widehat{S(t)v}(\xi)|^2d\xi &=\int_{\RR^d}(1+\norm{\xi}^2)^\beta E_{\alpha,1}^2(-t^\alpha\psi(\xi))|\widehat{v}(\xi)|^2d\xi\\
&\leq\int_{\RR^d}(1+\norm{\xi}^2)^\beta|\widehat{v}(\xi)|^2d\xi<\infty.
\end{align*}
Besides, denoting by $\langle\cdot,\cdot\rangle$ the usual Euclidean inner product on $\RR^d$, we obtain
\begin{align*}
AS(t)v&=\frac{1}{(2\pi)^d}\int_{\RR^d}e^{i\langle\cdot,\xi\rangle}(-\psi(\xi))\widehat{S(t)v}(\xi)d\xi\\
&=\frac{1}{(2\pi)^d}\int_{\RR^d}e^{i\langle\cdot,\xi\rangle}(-\psi(\xi))E_{\alpha,1}(-t^\alpha\psi(\xi))\widehat{v}(\xi)d\xi\\
&=\frac{1}{(2\pi)^d}\int_{\RR^d}e^{i\langle\cdot,\xi\rangle}E_{\alpha,1}(-t^\alpha\psi(\xi))(-\psi(\xi))\widehat{v}(\xi)d\xi\\
&=\frac{1}{(2\pi)^d}\int_{\RR^d}e^{i\langle\cdot,\xi\rangle}E_{\alpha,1}(-t^\alpha\psi(\xi))\widehat{Av}(\xi)d\xi\\
&=S(t)Av.
\end{align*}
The last condition $(S_3)$ also holds because
\begin{align*}
&\int_0^tg_\alpha(t-s)AS(s)vds\\
&=\int_0^tg_\alpha(t-s)\frac{1}{(2\pi)^d}\int_{\RR^d}e^{i\langle\cdot,\xi\rangle}(-\psi(\xi))E_{\alpha,1}(-s^\alpha\psi(\xi))\widehat{v}(\xi)d\xi ds\\
&=\frac{1}{(2\pi)^d}\int_{\RR^d}e^{i\langle\cdot,\xi\rangle}(-\psi(\xi))\widehat{v}(\xi)\int_0^t\frac{(t-s)^{\alpha-1}}{\Gamma(\alpha)}E_{\alpha,1}(-s^\alpha\psi(\xi))ds d\xi\\
&=\frac{1}{(2\pi)^d}\int_{\RR^d}e^{i\langle\cdot,\xi\rangle}(-\psi(\xi))\widehat{v}(\xi)\int_0^t\frac{(t-s)^{\alpha-1}}{\Gamma(\alpha)}\sum_{k=0}^\infty\frac{(-\psi(\xi))^k s^{\alpha k}}{\Gamma(\alpha k +1)} ds d\xi\\
&=\frac{1}{(2\pi)^d}\int_{\RR^d}e^{i\langle\cdot,\xi\rangle}(-\psi(\xi))\widehat{v}(\xi)\sum_{k=0}^\infty(-\psi(\xi))^k\int_0^t\frac{(t-s)^{\alpha-1} s^{\alpha k}}{\Gamma(\alpha)\Gamma(\alpha k +1)} ds d\xi\\
&=\frac{1}{(2\pi)^d}\int_{\RR^d}e^{i\langle\cdot,\xi\rangle}(-\psi(\xi))\widehat{v}(\xi)\sum_{k=0}^\infty(-\psi(\xi))^k\frac{t^{\alpha+\alpha k}}{\Gamma(\alpha k +1+\alpha)} d\xi\\
&=\frac{1}{(2\pi)^d}\int_{\RR^d}e^{i\langle\cdot,\xi\rangle}\widehat{v}(\xi)\sum_{k=0}^\infty\frac{(-\psi(\xi))^{k+1}(t^\alpha)^{k+1}}{\Gamma(\alpha(k +1)+1)} d\xi\\
&=\frac{1}{(2\pi)^d}\int_{\RR^d}e^{i\langle\cdot,\xi\rangle}\widehat{v}(\xi)\left[E_{\alpha,1}(-t^\alpha\psi(\xi))-1\right] d\xi\\
&=\frac{1}{(2\pi)^d}\int_{\RR^d}e^{i\langle\cdot,\xi\rangle}E_{\alpha,1}(-t^\alpha\psi(\xi))\widehat{v}(\xi)d\xi-\frac{1}{(2\pi)^d}\int_{\RR^d}e^{i\langle\cdot,\xi\rangle}\widehat{v}(\xi)d\xi\\
&=S(t)v-v.
\end{align*}
Since a strong solution $u$ of \eqref{Volterrageneral} is also a mild solution, from \cite[Proposition 1.2 part (i)]{Pru93} it follows that $u$ satisfies 
\begin{equation}
\label{formulaVP}
u(t)=\frac{d}{dt}\int_0^t S(s)f(t-s)ds, \quad t\in [0,T].
\end{equation}
The equation \eqref{formulaVP} is called the \textbf{variation of parameters formula} for the Volterra equation \eqref{Volterrageneral}.

Now, let $\widetilde{u}$ the unique local mild solution of \eqref{general} given by \eqref{Inteq:1} on $[0,T]$, under the assumptions of Theorem \ref{localsolution}. Define
\begin{equation}
\label{g}
g(t)(x)=g(t,x):=\lambda|\widetilde{u}(t,x)|^{\gamma-1}\widetilde{u}(t,x), \quad 0\leq t\leq T,\quad x\in\RR^d.
\end{equation}
We claim that $g_\alpha * g(t)\in L_2(\RR^d)$ for $0\leq t\leq T$, whenever $\frac{d\gamma}{\beta p}<1$. Indeed,
\begin{align*}
\norm{g_\alpha * g(t)}_{L_2(\RR^d)}&=\left(\int_{\RR^d}\left|(g_\alpha * g(t))(x)\right|^2dx\right)^{\frac{1}{2}}\\
&\leq\left(\int_{\RR^d}\left(\int_0^t\dfrac{(t-s)^{\alpha-1}}{\Gamma(\alpha)} |g(s,x)|ds\right)^2dx\right)^{\frac{1}{2}}\\
&\leq\int_0^t\dfrac{(t-s)^{\alpha-1}}{\Gamma(\alpha)}\left(\int_{\RR^d}|g(s,x)|^2dx\right)^{\frac{1}{2}}ds\\
&\lesssim\norm{\widetilde{u}}_{E_T}^{\gamma-1}\int_0^t(t-s)^{\alpha-1}s^{-\frac{\alpha d(\gamma-1)}{\beta p}}\left(\int_{\RR^d}|\widetilde{u}(s,x)|^2dx\right)^{\frac{1}{2}}ds\\
&\lesssim\norm{\widetilde{u}}_{E_T}^{\gamma-1}\int_0^t(t-s)^{\alpha-1}s^{-\frac{\alpha d(\gamma-1)}{\beta p}}\norm{\widetilde{u}(s,\cdot)}_2ds\\
&\lesssim\norm{\widetilde{u}}_{E_T}^{\gamma-1}\int_0^t(t-s)^{\alpha-1}s^{-\frac{\alpha d(\gamma-1)}{\beta p}}\max(\norm{\widetilde{u}(s,\cdot)}_1,\norm{\widetilde{u}(s,\cdot)}_\infty)ds\\
&\lesssim\norm{\widetilde{u}}_{E_T}^{\gamma-1}\int_0^t(t-s)^{\alpha-1}s^{-\frac{\alpha d(\gamma-1)}{\beta p}}(\norm{\widetilde{u}(s,\cdot)}_1+\norm{\widetilde{u}(s,\cdot)}_\infty)ds\\
&\lesssim\norm{\widetilde{u}}_{E_T}^{\gamma}\int_0^t(t-s)^{\alpha-1}s^{-\frac{\alpha d(\gamma-1)}{\beta p}}\left(1+s^{-\frac{\alpha d}{\beta p}}\right)ds\\
&\lesssim\norm{\widetilde{u}}_{E_T}^{\gamma}\left(t^{\alpha-\frac{\alpha d(\gamma-1)}{\beta p}}+t^{\alpha-\frac{\alpha d\gamma}{\beta p}}\right)\\
&<\infty.
\end{align*}
From here, we get
\[
g_\alpha * g(0)\equiv 0.
\]
In the same way it can be shown that $g_\alpha * g\in L_2([0,T];L_2(\RR^d))$. We additionally have that $g\in L_2([0,T];L_2(\RR^d))$ whenever $\frac{2\alpha d\gamma}{\beta p}<1$.

Let $u_0\in D(A)$. Consider the equations of Volterra
\begin{equation}
\label{EVg}
u(t)=g_\alpha * g(t)+\int_0^t g_\alpha(t-s)Au(s)ds
\end{equation}
and
\begin{equation}
\label{EVu0}
u(t)=u_0+\int_0^t g_\alpha(t-s)Au(s)ds,
\end{equation}
for $t\in[0,T]$.

Conclusions $(a)$ and $(b)$ of \cite[Theorem 8.7]{Pru93}, using $B\equiv 0$ and the Banach space $X_A=(H_2^\beta(\RR^d),\norm{\cdot}_{\Psi_{\beta},L_2})$, imply that \eqref{EVg} and \eqref{EVu0} have a unique a.e. strong solution $u_1$ and $u_2$, respectively.

\medbreak

Therefore, $u:=u_1+u_2\in L_2([0,T];X_A)$ is the unique a.e. strong solution of the Volterra equation
\begin{equation}
\label{ecuacionVolterra}
u(t)=u_0+g_\alpha * g(t)+\int_0^t g_\alpha(t-s)Au(s)ds
\end{equation} 
and satisfies \eqref{formulaVP} with $f(t):=u_0+g_\alpha * g(t)$, $t\in[0,T]$.

\medbreak

Next, we need to prove the following result using the operators $S(t)$ given by \eqref{Resolvente}.
\begin{lemma}
\label{PuntoFijoVP}
Let $\alpha\in (0,1)$ and $\beta\in (0,2)$. Assume the hypothesis $(\mathcal{H}_1)$ holds. Let $\lambda\in \RR$ and $\gamma>1$. Assume that $u$ is a fixed point of the operator given by \eqref{operator} on $E_T$ with $u_0\in L_2(\RR^d)$. If $\frac{d\gamma}{\beta p}<1$, then $u$ satisfies
\[
u(t)=\frac{d}{dt}\int_0^t S(s)f(t-s)ds,
\]
where $f(t):=u_0+g_\alpha * g(t)$ and $g$ like \eqref{g}, $t\in[0,T]$.
\end{lemma}
\begin{proof}
Let $x\in\RR^d$. We define
\begin{align*}
F(t,s)&:=\int_{\RR^d}g_{\alpha}*Z(\cdot,x-y)(t-s)g(s,y)dy\\
&=\int_{\RR^d}g(s,y)\left[\int_0^{t-s}\dfrac{(t-s-\tau)^{\alpha-1}}{\Gamma(\alpha)}Z(\tau,x-y)d\tau \right] dy\\
&=\int_0^{t-s}\dfrac{(t-s-\tau)^{\alpha-1}}{\Gamma(\alpha)}\int_{\RR^d}Z(\tau,x-y)g(s,y)dy d\tau\\
&\lesssim\norm{u}_{E_T}^\gamma\int_0^{t-s}\dfrac{(t-s-\tau)^{\alpha-1}}{\Gamma(\alpha)}s^{-\frac{\alpha d\gamma}{\beta p}}\underset{1}{\underbrace{\int_{\RR^d}Z(\tau,x-y)dy}} d\tau\\ 
&\lesssim\norm{u}_{E_T}^\gamma s^{-\frac{\alpha d\gamma}{\beta p}}(t-s)^\alpha.
\end{align*}
That is,
\begin{equation}
\label{condicion1}
g_{\alpha}*Z(\cdot,x-\cdot\cdot)(t-s)g(s,\cdot\cdot)\in L_1(\RR^d),\quad 0<s\leq t,\quad t>0.
\end{equation}
Without loss of generality, let $0<t_0<t$. The same arguments as in the proof of Theorem \ref{C1menosZ}, part $(i)$, but using the bounds of $Y$ given in Lemma \ref{cotasDeltaY}, yield
\begin{align*}
&\left|\int_{\RR^d} Y(t-s,x-y)g(s,y)dy-\int_{\RR^d} Y(t_0-s,x-y)g(s,y)dy\right|\\
&\leq \int_{\RR^d}\left|Y(t-s,x-y)-Y(t_0-s,x-y)\right||g(s,y)|dy\\
&\lesssim\norm{u}_{E_T}^\gamma s^{-\frac{\alpha d\gamma}{\beta p}}\int_{\RR^d}\left|Y(t-s,x-y)-Y(t_0-s,x-y)\right|dy\\
&\lesssim\norm{u}_{E_T}^\gamma s^{-\frac{\alpha d\gamma}{\beta p}}\norm{Y(t-s,x-\cdot)-Y(t_0-s,x-\cdot)}_1\\
&\lesssim\norm{u}_{E_T}^\gamma s^{-\frac{\alpha d\gamma}{\beta p}}|t-t_0|(t_0-s)^{\alpha-2}\rightarrow 0
\end{align*}
whenever $t\rightarrow t_0$.

This, together with \eqref{condicion1} and Lemma \ref{YenfunciondeZ}, proves that
\begin{equation}
\label{derivadaF}
\dfrac{\partial F}{\partial t}(t,s)=\int_{\RR^d} Y(t-s,x-y)g(s,y)dy,
\end{equation}
which is continuous w.r.t. $t\in(0, T]$, for $0<s<t$. Therefore, 
\[
\int_0^t\int_{\RR^d} Y(t-s,x-y)g(s,y)dyds=\int_0^t\dfrac{\partial F}{\partial t}(t,s)ds.
\]
From the previous work to \eqref{condicion1} we also obtain
\[
\left|\int_{t_0}^t F(t,s)ds\right|\lesssim\int_{t_0}^t s^{-\frac{\alpha d\gamma}{\beta p}}(t-s)^\alpha ds\lesssim t_0^{-\frac{\alpha d\gamma}{\beta p}}(t-t_0)^{\alpha+1}
\]
and thus
\begin{equation*}
\lim_{t\rightarrow t_0}\dfrac{1}{t-t_0}\int_{t_0}^t F(t,s)ds=0.
\end{equation*}
The mean-value theorem and the dominated convergence theorem yield
\begin{align*}
\frac{d}{dt}\left.\int_{0}^t F(t,s)ds\right|_{t=t_0}&=\lim_{t\rightarrow t_0}\dfrac{1}{t-t_0}\left(\int_{0}^t F(t,s)ds-\int_{0}^{t_0} F(t_0,s)ds\right)\\
&=\lim_{t\rightarrow t_0}\dfrac{1}{t-t_0}\int_{0}^{t_0}\left[F(t,s)-F(t_0,s)\right]ds+\lim_{t\rightarrow t_0}\dfrac{1}{t-t_0}\int_{t_0}^t F(t,s)ds\\
&=\lim_{t\rightarrow t_0}\int_{0}^{t_0}\frac{\partial F}{\partial t}(t_c,s)ds\\
&=\int_{0}^{t_0}\lim_{t\rightarrow t_0}\frac{\partial F}{\partial t}(t_c,s)ds\\
&=\int_{0}^{t_0}\frac{\partial F}{\partial t}(t_0,s)ds <\infty,
\end{align*}
because $t_0<t_c<t$ and the equality \eqref{derivadaF} implies
\begin{align*}
\left|\frac{\partial F}{\partial t}(t_c,s)\right|&\lesssim\int_{\RR^d}Y(t_c-s,x-y)|g(s,y)|dy\\
&\lesssim\norm{u}_{E_T}^\gamma s^{-\frac{\alpha d\gamma}{\beta p}}\int_{\RR^d}Y(t_c-s,x-y)dy\\
&\lesssim\norm{u}_{E_T}^\gamma s^{-\frac{\alpha d\gamma}{\beta p}}\frac{(t_c-s)^{\alpha-1}}{\Gamma(\alpha)}\\
&\lesssim\norm{u}_{E_T}^\gamma s^{-\frac{\alpha d\gamma}{\beta p}}(t_0-s)^{\alpha-1}.
\end{align*}
Consequently,
\begin{equation*}
\frac{d}{dt}\int_{0}^t F(t,s)ds=\int_{0}^{t}\frac{\partial F}{\partial t}(t,s)ds.
\end{equation*}
In view of \eqref{derivadaF}, the non-linear part of \eqref{operator} can be written as
\begin{align*}
\int_0^t\int_{\RR^d} Y(t-s,x-y)g(s,y)dyds&=\int_0^t\dfrac{\partial}{\partial t}F(t,s)ds\\
&=\frac{d}{dt}\int_{0}^t F(t,s)ds\\
&=\frac{d}{dt}\int_{0}^t \int_{\RR^d}g_{\alpha}*Z(\cdot,x-y)(t-s)g(s,y)dyds\\
&=\frac{d}{dt}\int_{0}^t \int_{\RR^d}g_{\alpha}*Z(\cdot,x-y)(s)g(t-s,y)dyds.
\end{align*}
Using Fubini and convolution respect to the time, the integral over $\RR^d$ is
\begin{equation*}
\int_0^s \int_{\RR^d} \dfrac{(s-\tau)^{\alpha-1}}{\Gamma(\alpha)}Z(\tau,x-y) g(t-s,y)dyd\tau
\end{equation*}
and hence
\begin{align*}
&\int_{0}^t \int_{\RR^d}g_{\alpha}*Z(\cdot,x-y)(s)g(t-s,y)dyds\\
&=\int_0^t\int_0^s \int_{\RR^d} \dfrac{(s-\tau)^{\alpha-1}}{\Gamma(\alpha)}Z(\tau,x-y) g(t-s,y)dyd\tau ds\\
&=\int_0^t\int_\tau^t \int_{\RR^d} \dfrac{(s-\tau)^{\alpha-1}}{\Gamma(\alpha)}Z(\tau,x-y) g(t-s,y)dydsd\tau\\
&=\int_0^t\int_{\RR^d}Z(\tau,x-y)\int_\tau^t \dfrac{(s-\tau)^{\alpha-1}}{\Gamma(\alpha)}g(t-s,y)dsdyd\tau\\
&=\int_0^t\int_{\RR^d}Z(\tau,x-y)\int_0^{t-\tau} \dfrac{(t-\tau-s)^{\alpha-1}}{\Gamma(\alpha)}g(s,y)dsdyd\tau\\
&=\int_0^t\int_{\RR^d}Z(\tau,x-y)g_\alpha * g(\cdot,y)(t-\tau)dyd\tau\\
&=\int_0^t Z(\tau,\cdot\cdot)\star (g_\alpha * g(\cdot,\cdot\cdot)(t-\tau))(x)d\tau\\
&=\int_0^t S(\tau)(g_\alpha * g(\cdot,\cdot\cdot)(t-\tau))(x)d\tau.
\end{align*}
Therefore, the operator \eqref{operator} evaluated in $u$ has the form
\[
u(t)=\frac{d}{dt}\int_0^t S(s)\left(u_0+g_\alpha * g(\cdot,\cdot\cdot)(t-s)\right)ds,\quad t\in[0,T].
\]
\end{proof}
By Lemma \ref{PuntoFijoVP} and uniqueness, we obtain $\widetilde{u}=u$. Besides, \eqref{ecuacionVolterra} is equivalent to
\begin{equation}
\label{localgeneral}
\begin{split}
\partial^{\alpha}_t(u-u_0)(t,x)+\Psi_{\beta}(-i\nabla)u(t,x) & = \lambda |u(t,x)|^{\gamma-1}u(t,x),\quad t\in(0,T], \; x\in\mathbb{R}^d,\\ u(t,x)|_{t=0} & = u_0(x),\quad x\in\mathbb{R}^d. 
\end{split}
\end{equation}
Indeed,   
\begin{align*}
&u(t)-\int_0^t g_\alpha(t-s)Au(s)ds=u_0+g_\alpha * g(t)\\
&\Leftrightarrow g_{1-\alpha}*u(t)+g_{1-\alpha}*g_\alpha*\Psi_{\beta}(-i\nabla)u(t)=g_{1-\alpha}*(u_0+g_\alpha * g)(t)\\
&\Leftrightarrow g_{1-\alpha}*(u-u_0)(t)+1*\Psi_{\beta}(-i\nabla)u(t)=1* g(t)\\
&\Leftrightarrow \frac{d}{dt}g_{1-\alpha}*(u-u_0)(t)+\Psi_{\beta}(-i\nabla)u(t)=g(t).
\end{align*}
Now, we define $u^+(t,x)=\max(u(t,x),0)$ and $u^-(t,x)=\max(-u(t,x),0)$.

From \cite[Section 2]{Zac08}, it is known that if $v\in L_2([0, T ];\RR)$, $g_{1-\alpha}*v\in W^1_2([0,T]; \RR)$ and $(g_{1-\alpha}*v)(0)=0$, then the operator $\partial^{\alpha}_t v:=\dfrac{d}{dt}(g_{1-\alpha}*v)$ has a Yosida approximation $\dfrac{d}{dt}(g_{1-\alpha,n}*v)$ in $L_2([0,T];\RR)$ as $n\rightarrow\infty$, with nonnegative and nonincreasing $g_{1-\alpha,n}\in W^1_1([0,T];\RR)$ for all $n\in\mathbb{N}$. From this work, one can also derive 
\begin{equation}
\label{CotaYosida}
v^-\dfrac{d}{dt}(g_{1-\alpha,n}*v)(t)\leq-\frac{1}{2}\frac{d}{dt}(g_{1-\alpha,n}*(v^-)^2)(t)\quad a.e.\; t\in(0,T),\quad n\in\mathbb{N}.
\end{equation}
Next, we use these results to prove that $u$ is a local positive solution a.e. of \eqref{general}, whenever $u_0\geq 0$, but non zero, and $\lambda<0$. By contradiction, suppose that $u<0$ somewhere on $(0,T]\times\RR^d$. Let $x\in\RR^d$. In order to apply \eqref{CotaYosida} to $u(\cdot,x)$, we need the condition $\frac{2\alpha d}{\beta p}<1$ to get $u(\cdot,x)\in L_2([0, T ];\RR)$. For the requirement $g_{1-\alpha}*u(\cdot,x)\in L_2([0,T]; \RR)$, we have that
\begin{align*}
\norm{g_{1-\alpha}*u(\cdot,x)}_2^2 &=\int_0^T|g_{1-\alpha}*u(t,x)|^2 dt\\
&=\int_0^T\left|\int_0^t \dfrac{(t-s)^{-\alpha}}{\Gamma(1-\alpha)}u(s,x)ds\right|^2 dt\\
&\lesssim\norm{u}_{E_T}^2\int_0^T\left[\int_0^t (t-s)^{-\alpha}s^{-\frac{\alpha d}{\beta p}}ds\right]^2 dt\lesssim\int_0^T\left[t^{1-\alpha-\frac{\alpha d}{\beta p}}\right]^2 dt
\end{align*}
which is finite because $0<\frac{2\alpha d}{\beta p}<1$ and $0<\alpha<1$.

We also find conditions such that $\dfrac{d}{dt}(g_{1-\alpha}*u)(\cdot,x)\in L_2([0, T ];\RR)$. For this purpose, we already know that $u\in L_2([0,T];X_A)$ is a strong solution of \eqref{ecuacionVolterra} and satisfies the equation \eqref{localgeneral}. Conclusions $(a)$ and $(b)$ of \cite[Theorem 8.7]{Pru93}) yield
\[
\int_0^T|\Psi_{\beta}(-i\nabla)u(t,x)|^2 dt<\infty\quad a.e.,\quad x\in\RR^d.
\] 
If $\frac{2\alpha d\gamma}{\beta p}<1$, we also get that $\abs{u}^{\gamma-1} u(\cdot,x)\in L_2([0, T ];\RR)$ and $\dfrac{d}{dt}(g_{1-\alpha}*u_0)(\cdot,x)\in L_2([0, T ];\RR)$ whenever $0<\alpha<\frac{1}{2}$. Finally, it can be readily checked that $(g_{1-\alpha}*u)(0,x)=0$ whenever $\alpha+\frac{\alpha d}{\beta p}<1$.

This work allows one to employ the Yosida approximation of $g_{1-\alpha}$ in $L_2([0, T ];\RR)$ with $u$. Using \eqref{CotaYosida} we obtain  
\[
u^-\dfrac{d}{dt}(g_{1-\alpha,n}*u)(t,x)\leq-\frac{1}{2}\frac{d}{dt}(g_{1-\alpha,n}*(u^-)^2)(t,x)\quad a.e.\; t\in(0,T).
\]
Besides, $(g_{1-\alpha,n}*(u^-)^2)(0,x)=0$ for all $n\in\mathbb{N}$ (see e.g., \cite[Formula 8]{Zac08} and \cite[Formula 10]{KSVZ16}). Consequently, 
\begin{equation*}
\begin{split}
&\int_0^T u^-\dfrac{d}{dt}(g_{1-\alpha,n}*u)(t,x)dt\\
&\leq-\frac{1}{2}\int_0^T\frac{d}{dt}(g_{1-\alpha,n}*(u^-)^2)(t,x)dt\\
&=-\frac{1}{2}(g_{1-\alpha,n}*(u^-)^2)(T,x)+\frac{1}{2}(g_{1-\alpha,n}*(u^-)^2)(0,x)\\
&=-\frac{1}{2}(g_{1-\alpha,n}*(u^-)^2)(T,x)\\
&\leq 0.
\end{split}
\end{equation*}
Thus,
\begin{equation*}
\begin{split}
\int_0^T u^-\dfrac{d}{dt}(g_{1-\alpha}*u)(t,x)dt \leq\int_0^T u^-\frac{d}{dt}((g_{1-\alpha}-g_{1-\alpha,n})*u)(t,x)dt
\end{split}
\end{equation*}
and applying Hölder we conclude that
\begin{equation*}
\begin{split}
&\int_0^T u^-\dfrac{d}{dt}(g_{1-\alpha}*u)(t,x)dt\\
& \leq\left(\int_0^T (u^-)^2(t,x) dt\right)^{\frac{1}{2}}\left(\int_0^T\left|\frac{d}{dt}((g_{1-\alpha}-g_{1-\alpha,n})*u)(t,x)\right|^2 dt\right)^{\frac{1}{2}}\\
& < \norm{u^-(\cdot,x)}_2 \varepsilon.
\end{split}
\end{equation*}
This shows that 
\begin{equation}
\label{DerivadaFraccionalNegativa}
\int_0^T u^-\dfrac{d}{dt}(g_{1-\alpha}*(u-u_0))(t,x)dt \leq 0
\end{equation}
a.e., $x\in\RR^d$.

On the other hand, if $v, w\in E_T$ we have that $v(t,\cdot), w(t,\cdot)\in L_2(\RR^d)$ whenever $t>0$. Due to Courrège, $-\Psi_{\beta}(-i\nabla)$ satisfies \cite[Theorem 3.6]{Schi01} and we can write
\begin{equation}
\label{ProductoInterno}
\langle \Psi_{\beta}(-i\nabla) v(t,\cdot),w(t,\cdot)\rangle_2=\frac{1}{2}\int_{\RR^{2d}}(v(t,y)-v(t,x))(w(t,y)-w(t,x))\nu(dy)dx
\end{equation}
with a positive kernel $\nu$ on $\RR^d$ such that $\displaystyle\int_{\RR^d\setminus\{x\}}\frac{\norm{y-x}^2}{1+\norm{y-x}^2}\nu(dy)<\infty$. Here, we use the notation $\langle \cdot, \cdot\rangle_2$ to refer to the inner product on $L_2(\RR^d)$, i.e.,  
\[
\langle \vartheta_1,\vartheta_2\rangle_2=\int_{\RR^d}\vartheta_1(x)\vartheta_2(x)dx,\quad \vartheta_1,\vartheta_2\in L_2(\RR^d).
\]
Using \eqref{localgeneral}, we obtain
\begin{equation*}
\begin{split}
u^-(t,x)\partial^{\alpha}_t (u-u_0)(t,x)=&-u^-(t,x)\Psi_{\beta}(-i\nabla)u(t,x) + \lambda u^-(t,x)\abs{u(t,x)}^{\gamma-1} u(t,x)\\
=&-u^-(t,x)\Psi_{\beta}(-i\nabla)u^+(t,x)+u^-(t,x)\Psi_{\beta}(-i\nabla)u^-(t,x)\\
&- \lambda  \abs{u(t,x)}^{\gamma-1} (u^-)^2(t,x).
\end{split}
\end{equation*}
Integrating over $[0,T]\times\RR^d$ and using Fubini in a convenient form, we get
\begin{equation*}
\begin{split}
\int_{\RR^d}\int_0^T u^-(t,x)\partial^{\alpha}_t (u-u_0)(t,x)dt dx=&-\int_0^T\int_{\RR^d} u^-(t,x)\Psi_{\beta}(-i\nabla)u^+(t,x)dx dt\\
&+\int_0^T\int_{\RR^d}u^-(t,x)\Psi_{\beta}(-i\nabla)u^-(t,x)dx dt\\
&-\lambda \int_0^T\int_{\RR^d} \abs{u(t,x)}^{\gamma-1} (u^-)^2(t,x)dx dt.
\end{split}
\end{equation*}
From \eqref{DerivadaFraccionalNegativa} and \eqref{ProductoInterno} it follows that the left-hand side is non-positive and the right-hand side is strictly positive, respectively. This contradiction shows that $u\geq 0$ a.e. For the case $\frac{1}{2}\leq \alpha <1$ we set the parameter $\bar{p}$, such that $1<\bar{p}<\frac{1}{\alpha}$. Due to this choice, we have that $\bar{p}<\frac{\bar{p}}{\bar{p}-1}$ and that $2<\frac{\bar{p}}{\bar{p}-1}$. Thereby, we can apply similar arguments, i.e., Yosida approximation of $g_{1-\alpha}$ in $L_{\bar{p}}([0, T ];\RR)$ and \cite[Theorem 8.7]{Pru93} in $L_{\bar{p}}([0, T ];L_2(\RR^d))$, for obtaining again that $u\geq 0$ a.e. Whenever $0<\alpha <\frac{1}{2}$, we fix $\bar{p}=2$. We are now in a position to show the following theorem.
\begin{theorem}
\label{globalpositivesolution}
Let $\alpha\in (0,1)$ and $\beta\in (0,2)$. Assume the hypothesis $(\mathcal{H}_1)$ holds. Let $\lambda < 0$, $\gamma>1$. Suppose that $\max\left(1,\kappa,\frac{d\gamma}{\beta}\right)<p<\infty$, that $\frac{\bar{p}\alpha d}{(\bar{p}-1)\beta p}<1$, that $\frac{\bar{p}\alpha d\gamma}{\beta p}<1$ and that $\alpha+\frac{\alpha d}{\beta p}<1$. If $u_0\in L_p(\RR^d) \cap L_1(\RR^d)\cap H_2^\beta(\RR^d)$ is non-negative a.e., then there exists a unique non-negative global solution $u\in C([0,\infty); L_p(\RR^d) \cap L_1(\RR^d)) \cap C((0,\infty); L_{\infty}(\RR^d))$ to the Cauchy problem \eqref{general}. Moreover, estimate
\[
\norm{u(t)}_1+t^{\frac{\alpha d}{\beta}\left(\frac{1}{p'}-\frac{1}{p}\right)}\norm{u(t)}_p+t^{\frac{\alpha d}{\beta p'}}\norm{u(t)}_\infty\lesssim(\norm{u_0}_1+\norm{u_0}_p) 
\]  
is true for all $t\geq 1$, with $p'$ as in Theorem \ref{globalsolutionsmallu0}. 
\end{theorem}

\begin{proof}
As in Sect.~\ref{sec:2}, one finds that there exists a unique local solution $u\in E_{T^*}$ for some $T^*>0$. This solution is non-negative a.e., as discussed above, and satisfies
\[
u(t,x)=\int_{\RR^d} Z(t,x-y)u_0(y)dy + \lambda\int_0^t\int_{\RR^d} Y(t-s,x-y)|u(s,y)|^{\gamma-1}u(s,y)dyds.
\]
The positivity of $Z$ and $Y$ yields
\[
u(t,x)\leq\int_{\RR^d} Z(t,x-y)u_0(y)dy.
\]
Therefore,

\medbreak

$\norm{u(t,\cdot)}_1\leq\left\|Z(t,\cdot)\star u_0\right\|_1\leq\norm{u_0}_1$ for $t\in [0,T^*]$,

\medbreak

$\norm{u(t,\cdot)}_p\leq\left\|Z(t,\cdot)\star u_0\right\|_p\leq\norm{u_0}_p$ for $t\in [0,T^*]$,

\medbreak

$\norm{u(t,\cdot)}_\infty\leq\left\|Z(t,\cdot)\star u_0\right\|_\infty\lesssim t^{-\frac{\alpha d}{\beta p}}\norm{u_0}_p$ for $t\in (0,T^*]$

\medbreak

and thus we obtain
\begin{equation}
\label{cotaglobal}
\norm{u}_{E_{T^*}}\leq C \norm{u_0}_{L_p(\RR^d) \cap L_1(\RR^d)}.
\end{equation}
We note that the constant $C$ is independent of $T^*$, that is, we can apply Theorem \ref{localsolution} on $[T^*,T_1]$, with $0<T^*<T_1$, and the extended solution $u\in E_{T_1}$ also satisfies \eqref{cotaglobal}.

In particular, $\norm{u(T^*)}_{L_p(\RR^d) \cap L_1(\RR^d)}\leq C \norm{u_0}_{L_p(\RR^d) \cap L_1(\RR^d)}$.  This estimate allows to prolong the local solution for all times $t>0$ (see \cite[Theorem 1.20]{HKNS06}). Since $u(T^*)$ is non-negative and it is the new initial condition on $[T^*,T_1]$, it follows that $u\in E_{T_1}$ is also non-negative. Consequently, the global solution is non-negative. 

On the other hand, as in proof of Theorem \ref{globalsolutionsmallu0}, we also get for all $t\geq 1$,
\begin{align*}
\norm{u(t,\cdot)}_p\leq\left\|Z(t,\cdot)\star u_0\right\|_p\lesssim t^{-\frac{\alpha d}{\beta}\left(\frac{1}{p'}-\frac{1}{p}\right)}\max (\norm{u_0}_1,\norm{u_0}_p),\\
\norm{u(t,\cdot)}_\infty\leq\left\|Z(t,\cdot)\star u_0\right\|_{\infty}\lesssim t^{-\frac{\alpha d}{\beta p'}}\max (\norm{u_0}_1,\norm{u_0}_p).
\end{align*}
Thus, we obtain the desired estimate.
\end{proof}

\section{Asymptotic behavior}
\label{sec:5}
In this section we study the $L_p$-decay of a global solution, which was obtained in Sect.~\ref{sec:3} and in Sect.~\ref{sec:4}, respectively. We recall that in the first case (Theorem \ref{globalsolutionsmallu0}), a small initial data is required. In the other case (Theorem \ref{globalpositivesolution}), $\lambda<0$ and initial data non-negative are required, which yield a global solution non-negative.
\begin{lemma}
\label{NormaGradienteZ}
Let $\alpha\in (0,1)$ and $\beta\in (1,2)$. Assume the hypothesis $(\mathcal{H}_2)$ holds. Then there exists a positive constant $C$ for all $t>0$ and $y\in\RR^d$, such that the estimate
\[		
\norm{Z(t,\cdot-y)-Z(t,\cdot)}_q\leq C \norm{y}\norm{\nabla Z(t,\cdot)}_q\lesssim \norm{y} t^{-\frac{\alpha d}{\beta}\left(1-\frac{1}{q}\right)-\frac{\alpha}{\beta}}	
\]
is true for $1\leq q<\frac{d}{d+1-\beta}$.	
\end{lemma}

\begin{proof}
By application of bounds given in Lemma \ref{cotasDeltaEspacialZ}, we have
	\[
		|Z(t,x-y)-Z(t,x)|\leq C \norm{y}\begin{cases}t^{-\alpha}\norm{x-\varepsilon y}^{\beta-(d+1)}, &\norm{x-\varepsilon y}\leq t^{\frac{\alpha}{\beta}} \\t^{\alpha}\norm{x-\varepsilon y}^{-\beta-(d+1)}, &\norm{x-\varepsilon y}\geq t^{\frac{\alpha}{\beta}} 		
		\end{cases}=:C \norm{y}D(t,x-\varepsilon y),
\]
where $\varepsilon\in (0,1)$.

Using an argument as in the Hardy's inequality proof, we get
\begin{align*}
\norm{Z(t,\cdot-y)-Z(t,\cdot)}_q &=\left(\int_{\RR^d}|Z(t,x-y)-Z(t,x)|^q dx\right)^{\frac{1}{q}}\\
&\lesssim\norm{y} \left(\int_{\RR^d}D^q(t,x-\varepsilon y) dx\right)^{\frac{1}{q}}\\
&\lesssim\norm{y} \left(\int_{\RR^d}\left(\int_0^1 D(t,x-s y)ds\right)^q dx\right)^{\frac{1}{q}}\\
&\lesssim\norm{y} \int_0^1\left(\int_{\RR^d} D^q(t,x-s y) dx\right)^{\frac{1}{q}}ds\\
&\lesssim\norm{y}\norm{D(t,\cdot)}_q.
\end{align*}
Now, $\norm{D(t,\cdot)}_q$ can be estimated in the same way as $Z(t,\cdot)$ in Theorem \ref{CotaNormapZ} and $\norm{\nabla Z(t,\cdot)}_q\lesssim \norm{D(t,\cdot)}_q$. The mean-value inequality completes the proof.
\end{proof}
On the other hand, from \cite{DZ92} is known the following decomposition lemma.
\begin{lemma}
Assume that $1\leq r < \frac{d}{d-1}$, $f\in L_1(\RR^d)$ and $\norm{\cdot}f\in L_r(\RR^d)$, then there exists
a vectorial function $\textbf{F}\in L_r(\RR^d;\RR^d)$ such that 
\[
f=\left(\int_{\RR^d}f(y)dy\right)\delta_0+div\textbf{F}.
\]
in the distributional sense and
\[
\norm{\textbf{F}}_r\leq C(d)\norm{\norm{\cdot}f}_r.
\]
\end{lemma}
Using this with $f=u_0$, we find that
\begin{align*}
Z(t,\cdot)\star u_0(x)&=\left(\int_{\RR^d}u_0(y)dy\right)Z(t,\cdot)\star\delta_0(x)+Z(t,\cdot)\star div\textit{\textbf{F}}(x)\\
&=\left(\int_{\RR^d}u_0(y)dy\right)Z(t,x)+\nabla Z(t,\cdot)\star \textit{\textbf{F}}(x).
\end{align*}
Let $A=\displaystyle\int_{\RR^d} u_0(y)dy$. Young’s inequality for convolutions yields
\begin{align*}
\left\|Z(t,\cdot)\star u_0-A Z(t,\cdot)\right\|_p&\lesssim \norm{\nabla Z(t,\cdot)}_{q}\norm{\textit{\textbf{F}}}_r\\
&\lesssim \norm{\nabla Z(t,\cdot)}_{q}\norm{\norm{\cdot}u_0}_r
\end{align*}
and we get
\begin{equation}
\label{decayinit}
\norm{Z(t,\cdot)\star u_0-AZ(t,\cdot)}_p\lesssim t^{-\frac{\alpha d}{\beta}\left(1-\frac{1}{q}\right)-\frac{\alpha}{\beta}}\norm{\norm{\cdot}u_0}_r.
\end{equation}
Now, the idea is to find conditions such that
\begin{equation}
\label{pqr}
\frac{1}{p}+1=\frac{1}{q}+\frac{1}{r},\quad 1\leq q<\frac{d}{d+1-\beta},\quad 1\leq r<\frac{d}{d-1}.
\end{equation}
First, we need $\beta\in (1,2)$ since Lemma \ref{NormaGradienteZ} requires this for choices of $q$. Besides, in the case $d=1$ is enough $p<\infty$ but in the case $d\geq 2$ we need $p<\frac{d}{d-\beta}$. In this way it is possible to get the estimate \eqref{decayinit} and the required upper bound on $p$ it allows one to achieve both the $L_p$-norm and the $L_q$-norm for $Y$.

Therefore, the estimate 
\begin{equation}
\label{decayforc}
\left\|\lambda\int_0^t Y(t-s,\cdot\cdot)\star|u(s,\cdot\cdot)|^{\gamma-1}u(s,\cdot\cdot)(\cdot)ds-BY(t,\cdot)\right\|_p\lesssim t^{-\frac{\alpha d}{\beta}\left(\frac{1}{p'}-\frac{1}{p}\right)+\alpha-1}
\end{equation}
is true for $t\geq 1$, with $B=\lambda\displaystyle\int_0^\infty\int_{\RR^d}|u(s,y)|^{\gamma-1}u(s,y)dyds$, whenever $\frac{\alpha d}{\beta p'}(\gamma-1)>1$. Here, $p'$ is the same as in Theorems \ref{globalsolutionsmallu0} and \ref{globalpositivesolution}. Indeed,
\begin{align*}
&\int_0^t Y(t-s,\cdot\cdot)\star|u(s,\cdot\cdot)|^{\gamma-1}u(s,\cdot\cdot)(x)ds-Y(t,x)\int_0^\infty\int_{\RR^d}|u(s,y)|^{\gamma-1}u(s,y)dyds\\
=&\int_0^{\frac{t}{2}} Y(t-s,\cdot\cdot)\star|u(s,\cdot\cdot)|^{\gamma-1}u(s,\cdot\cdot)(x)ds + \int_{\frac{t}{2}}^t Y(t-s,\cdot\cdot)\star|u(s,\cdot\cdot)|^{\gamma-1}u(s,\cdot\cdot)(x)ds\\
&-\displaystyle\int_0^{\frac{t}{2}}(Y(t,x)-Y(t-s,x)+Y(t-s,x))\int_{\RR^d}|u(s,y)|^{\gamma-1}u(s,y)dyds\\
&-\displaystyle\int_{\frac{t}{2}}^\infty Y(t,x)\int_{\RR^d}|u(s,y)|^{\gamma-1}u(s,y)dyds\\
=&\int_0^{\frac{t}{2}} Y(t-s,\cdot\cdot)\star|u(s,\cdot\cdot)|^{\gamma-1}u(s,\cdot\cdot)(x)ds + \int_{\frac{t}{2}}^t Y(t-s,\cdot\cdot)\star|u(s,\cdot\cdot)|^{\gamma-1}u(s,\cdot\cdot)(x)ds\\
&-\displaystyle\int_0^{\frac{t}{2}}Y(t-s,x)\int_{\RR^d}|u(s,y)|^{\gamma-1}u(s,y)dyds-\displaystyle\int_{\frac{t}{2}}^\infty Y(t,x)\int_{\RR^d}|u(s,y)|^{\gamma-1}u(s,y)dyds\\
&-\displaystyle\int_0^{\frac{t}{2}}(Y(t,x)-Y(t-s,x))\int_{\RR^d}|u(s,y)|^{\gamma-1}u(s,y)dyds\\
=:&J_1(t,x)+J_2(t,x)-J_3(t,x)-J_4(t,x)-J_5(t,x).
\end{align*}
Now, we estimate the $L_p$-norm for $J_k$, $k=1,2,3,4,5$, recalling that $\frac{d}{\beta p}(\gamma-1)<1$, $\frac{d}{\beta p'}<1$ and $1\leq p'<p$. 
\begin{align*}
\norm{J_1(t,\cdot)}_p & \leq \int_0^{\frac{t}{2}} \norm{Y(t-s,\cdot\cdot)\star|u(s,\cdot\cdot)|^{\gamma-1}u(s,\cdot\cdot)}_p ds\\
& \leq \int_0^{\frac{t}{2}} \norm{Y(t-s,\cdot\cdot)}_p\norm{|u(s,\cdot\cdot)|^{\gamma-1}u(s,\cdot\cdot)}_1 ds\\
& \lesssim\norm{u}_{E}^\gamma \int_0^{\frac{t}{2}} (t-s)^{-\frac{\alpha d}{\beta}\left(1-\frac{1}{p}\right)+\alpha-1}s^{-\frac{\alpha d}{\beta p}(\gamma-1)}\langle s\rangle^{-\frac{\alpha d}{\beta}(\gamma-1)\left(\frac{1}{p'}-\frac{1}{p}\right)} ds\\
& \lesssim\norm{u}_{E}^\gamma\left(\frac{t}{2}\right)^{-\frac{\alpha d}{\beta}\left(1-\frac{1}{p}\right)+\alpha-1} \int_0^{\frac{t}{2}} s^{-\frac{\alpha d}{\beta p}(\gamma-1)}\langle s\rangle^{-\frac{\alpha d}{\beta}(\gamma-1)\left(\frac{1}{p'}-\frac{1}{p}\right)} ds.
\end{align*}
Choosing $0<c<\frac{t}{2}$ we get
\begin{align*}
\int_0^{\frac{t}{2}} s^{-\frac{\alpha d}{\beta p}(\gamma-1)}\langle s\rangle^{-\frac{\alpha d}{\beta}(\gamma-1)\left(\frac{1}{p'}-\frac{1}{p}\right)} ds&=\int_0^c s^{-\frac{\alpha d}{\beta p}(\gamma-1)}\langle s\rangle^{-\frac{\alpha d}{\beta}(\gamma-1)\left(\frac{1}{p'}-\frac{1}{p}\right)} ds+\int_c^{\frac{t}{2}}s^{-\frac{\alpha d}{\beta p}(\gamma-1)}\langle s\rangle^{-\frac{\alpha d}{\beta}(\gamma-1)\left(\frac{1}{p'}-\frac{1}{p}\right)} ds\\
&\lesssim\int_0^c s^{-\frac{\alpha d}{\beta p}(\gamma-1)}ds+\int_c^{\frac{t}{2}}s^{-\frac{\alpha d}{\beta p'}(\gamma-1)} ds\\
&\lesssim\frac{c^{1-\frac{\alpha d}{\beta p}(\gamma-1)}}{1-\frac{\alpha d}{\beta p}(\gamma-1)}+\frac{1}{1-\frac{\alpha d}{\beta p'}(\gamma-1)}\left(\left(\frac{t}{2}\right)^{1-\frac{\alpha d}{\beta p'}(\gamma-1)}-c^{1-\frac{\alpha d}{\beta p'}(\gamma-1)}\right).
\end{align*}
We note that $1-\frac{\alpha d}{\beta p}(\gamma-1)>0$ and $1-\frac{\alpha d}{\beta p'}(\gamma-1)<0$. Thus,
\begin{align*}
\norm{J_1(t,\cdot)}_p &\lesssim \norm{u}_{E}^\gamma\left(\frac{t}{2}\right)^{-\frac{\alpha d}{\beta}\left(1-\frac{1}{p}\right)+\alpha-1}\\
&\lesssim t^{-\frac{\alpha d}{\beta}\left(1-\frac{1}{p}\right)+\alpha-1}.
\end{align*}
However, if $t\geq 1$ we have that $t^{-\frac{\alpha d}{\beta}\left(1-\frac{1}{p}\right)}\leq t^{-\frac{\alpha d}{\beta}\left(\frac{1}{p'}-\frac{1}{p}\right)}$ and we obtain
\begin{equation}
\label{J_1}
\norm{J_1(t,\cdot)}_p\lesssim t^{-\frac{\alpha d}{\beta}\left(\frac{1}{p'}-\frac{1}{p}\right)+\alpha-1}.
\end{equation}
We continue with $J_2$ and we estimate
\begin{align*}
\norm{J_2(t,\cdot)}_p & \leq \int_{\frac{t}{2}}^t \norm{Y(t-s,\cdot\cdot)\star|u(s,\cdot\cdot)|^{\gamma-1}u(s,\cdot\cdot)}_p ds\\
& \leq \int_{\frac{t}{2}}^t \norm{Y(t-s,\cdot\cdot)}_q\norm{|u(s,\cdot\cdot)|^{\gamma-1}u(s,\cdot\cdot)}_{p'} ds\\
& \lesssim\norm{u}_{E}^\gamma \int_{\frac{t}{2}}^t (t-s)^{-\frac{\alpha d}{\beta}\left(\frac{1}{p'}-\frac{1}{p}\right)+\alpha-1}s^{-\frac{\alpha d}{\beta p'}(\gamma-1)}ds\\
& \lesssim\norm{u}_{E}^\gamma\left(\frac{t}{2}\right)^{-\frac{\alpha d}{\beta p'}(\gamma-1)} \int_{\frac{t}{2}}^t (t-s)^{-\frac{\alpha d}{\beta}\left(\frac{1}{p'}-\frac{1}{p}\right)+\alpha-1}ds\\
& \lesssim\norm{u}_{E}^\gamma\left(\frac{t}{2}\right)^{-\frac{\alpha d}{\beta p'}(\gamma-1)} \left. (t-s)^{-\frac{\alpha d}{\beta}\left(\frac{1}{p'}-\frac{1}{p}\right)+\alpha}\right|_{\frac{t}{2}}^t\\
&\lesssim\norm{u}_{E}^\gamma\left(\frac{t}{2}\right)^{-\frac{\alpha d}{\beta p'}(\gamma-1)} \left(\frac{t}{2}\right)^{-\frac{\alpha d}{\beta}\left(\frac{1}{p'}-\frac{1}{p}\right)+\alpha}.
\end{align*}
Therefore,
\begin{equation}
\label{J_2}
\norm{J_2(t,\cdot)}_p\lesssim t^{-\frac{\alpha d}{\beta}\left(\frac{1}{p'}-\frac{1}{p}\right)+\alpha-\frac{\alpha d}{\beta p'}(\gamma-1)}.
\end{equation} 
For $J_3$ we have that
\begin{align*}
\norm{J_3(t,\cdot)}_p & \leq \displaystyle\int_0^{\frac{t}{2}}\norm{Y(t-s,\cdot)}_p\int_{\RR^d}|u(s,y)|^{\gamma-1}|u(s,y)|dy\;ds\\
&\lesssim\norm{u}_{E}^\gamma\int_0^{\frac{t}{2}} (t-s)^{-\frac{\alpha d}{\beta}\left(1-\frac{1}{p}\right)+\alpha-1}s^{-\frac{\alpha d}{\beta p}(\gamma-1)}\langle s\rangle^{-\frac{\alpha d}{\beta}(\gamma-1)\left(\frac{1}{p'}-\frac{1}{p}\right)} ds\\
& \lesssim\norm{u}_{E}^\gamma\left(\frac{t}{2}\right)^{-\frac{\alpha d}{\beta}\left(1-\frac{1}{p}\right)+\alpha-1} \int_0^{\frac{t}{2}} s^{-\frac{\alpha d}{\beta p}(\gamma-1)}\langle s\rangle^{-\frac{\alpha d}{\beta}(\gamma-1)\left(\frac{1}{p'}-\frac{1}{p}\right)} ds
\end{align*}
and proceeding in the same way as in the estimate of $J_1$, we obtain
\begin{equation}
\label{J_3}
\norm{J_3(t,\cdot)}_p\lesssim t^{-\frac{\alpha d}{\beta}\left(\frac{1}{p'}-\frac{1}{p}\right)+\alpha-1}.
\end{equation} 
For $J_4$ we find
\begin{align*}
\norm{J_4(t,\cdot)}_p & \leq \norm{Y(t,\cdot)}_p\displaystyle\int_{\frac{t}{2}}^\infty\int_{\RR^d}|u(s,y)|^{\gamma-1}|u(s,y)|dyds\\
&\lesssim\norm{u}_{E}^\gamma t^{-\frac{\alpha d}{\beta}\left(1-\frac{1}{p}\right)+\alpha-1}\int_{\frac{t}{2}}^\infty s^{-\frac{\alpha d}{\beta p'}(\gamma-1)} ds\\
&\lesssim\norm{u}_{E}^\gamma t^{-\frac{\alpha d}{\beta}\left(1-\frac{1}{p}\right)+\alpha-1}\left(\frac{t}{2}\right)^{1-\frac{\alpha d}{\beta p'}(\gamma-1)}
\end{align*}
and we get
\begin{equation}
\label{J_4}
\norm{J_4(t,\cdot)}_p\lesssim t^{-\frac{\alpha d}{\beta}\left(\frac{1}{p'}-\frac{1}{p}\right)+\alpha-\frac{\alpha d}{\beta p'}(\gamma-1)}.
\end{equation} 
In order to estimate $J_5$, we follow the same arguments as in the proof of Theorem \ref{C1menosZ}, first part, but using the bounds of $Y$ given in Lemma \ref{cotasDeltaY}, that is,
\begin{align*}
\norm{J_5(t,\cdot)}_p & \leq \displaystyle\int_0^{\frac{t}{2}}\norm{Y(t,\cdot)-Y(t-s,\cdot)}_p\int_{\RR^d}|u(s,y)|^{\gamma-1}|u(s,y)|dyds\\
& \lesssim\norm{u}_{E}^\gamma\int_0^{\frac{t}{2}}s (t-s)^{-\frac{\alpha d}{\beta}\left(1-\frac{1}{p}\right)+\alpha-2}s^{-\frac{\alpha d}{\beta p}(\gamma-1)}\langle s\rangle^{-\frac{\alpha d}{\beta}(\gamma-1)\left(\frac{1}{p'}-\frac{1}{p}\right)} ds\\
& \lesssim\norm{u}_{E}^\gamma\left(\frac{t}{2}\right)^{-\frac{\alpha d}{\beta}\left(1-\frac{1}{p}\right)+\alpha-1} \int_0^{\frac{t}{2}} s(t-s)^{-1}s^{-\frac{\alpha d}{\beta p}(\gamma-1)}\langle s\rangle^{-\frac{\alpha d}{\beta}(\gamma-1)\left(\frac{1}{p'}-\frac{1}{p}\right)} ds\\
& \lesssim\norm{u}_{E}^\gamma\left(\frac{t}{2}\right)^{-\frac{\alpha d}{\beta}\left(1-\frac{1}{p}\right)+\alpha-1} \int_0^{\frac{t}{2}} s\;s^{-1}s^{-\frac{\alpha d}{\beta p}(\gamma-1)}\langle s\rangle^{-\frac{\alpha d}{\beta}(\gamma-1)\left(\frac{1}{p'}-\frac{1}{p}\right)} ds\\
&\lesssim\norm{u}_{E}^\gamma\left(\frac{t}{2}\right)^{-\frac{\alpha d}{\beta}\left(1-\frac{1}{p}\right)+\alpha-1} \int_0^{\frac{t}{2}} s^{-\frac{\alpha d}{\beta p}(\gamma-1)}\langle s\rangle^{-\frac{\alpha d}{\beta}(\gamma-1)\left(\frac{1}{p'}-\frac{1}{p}\right)} ds.
\end{align*}
Now, we proceed in the same way as in the estimate of $J_1$ and thus
\begin{equation}
\label{J_5}
\norm{J_5(t,\cdot)}_p\lesssim t^{-\frac{\alpha d}{\beta}\left(\frac{1}{p'}-\frac{1}{p}\right)+\alpha-1}.
\end{equation} 
Gathering estimates from \eqref{J_1} to \eqref{J_5}, we have proved \eqref{decayforc}. This, together with \eqref{decayinit}, show the following result.   
\begin{theorem}
\label{decay}
Let $\alpha\in (0,1)$ and $\beta\in (1,2)$. Assume the hypothesis $(\mathcal{H}_2)$ holds. Let $\lambda\in \RR$ and $\gamma>1$. Suppose that the initial data $u_0\in L_1(\RR^d)$, $p$  and $p'$ satisfy conditions as in Theorem \ref{globalsolutionsmallu0}. Let $p<\kappa_1$ as in Theorem \ref{CotaNormapZ} and $\frac{\alpha d}{\beta p'}(\gamma-1)>1$. Assume in addition that $\norm{\cdot}u_0\in L_r(\RR^d)$, with some $r$ and $q$ as in \eqref{pqr}. If $u\in E$ is a global solution to the Cauchy problem \eqref{general}, then $u$ has the asymptotic behavior
\[
\norm{u(t,\cdot)-AZ(t,\cdot)-BY(t,\cdot)}_p\rightarrow 0
\]
as $t\rightarrow\infty$, with the constants
\[
A=\displaystyle\int_{\RR^d} u_0(y)dy
\]
and
\[
B=\lambda\displaystyle\int_0^\infty\int_{\RR^d}|u(s,y)|^{\gamma-1}u(s,y)dyds.
\]
\end{theorem}


%
\section*{Funding}
The authors were partially supported by Chilean research grant Fondo Nacional de Desarrollo Científico y Tecnológico, FONDECYT 1190255.

\section*{Conflict of interest}
The authors declare that they have no conflict of interest.

\section*{Availability of data and material}
Not applicable. No datasets were generated or analysed during the current study.

\section*{Code availability}
Not applicable.

\bibliographystyle{ieeetr}

\bibliography{bibl}

\begin{thebibliography}{99}

\bibitem {Aro67} D. Aronson. \textit{Bounds for the fundamental solution of a parabolic equation}.  Bull. Amer. Math. Soc.,  Vol. 73, No. 6: 890-896, 1967.

\bibitem {BE89} J. Bebernes, D. Eberly. \textit{Mathematical problems from combustion theory}. Applied Mathematical Sciences, vol 83. Springer, New York, 1989.

\bibitem {Bru58} N. G. De Bruijin. \textit{Asymptotic methods in analysis}. North-Holland Publishing Co. - Amsterdam P. Noordhoff Ltd. - Groningen, Netherlands, 1958.

\bibitem {Dav07} E. Davies. \textit{Linear operators and their spectra}. Cambridge University Press, 2007.

\bibitem {DZ92} J. Duoandikoetxea, E. Zuazua. \textit{Moments, masses de Diracet décomposition de fonctions}. C. R. Acad. Sci. Sér. 1 Math. 315(6):693–698, 1992.
 
\bibitem {EIK04} S.D. Eidelman, S.D. Ivasyshen, A.N. Kochubei. \textit{Analytic Methods in the Theory of Differential and Pseudo-Differential Equations of Parabolic Type}. Springer Basel AG, 2004.

\bibitem {EK04} S.D. Eidelman, A.N. Kochubei. \textit{Cauchy problem for fractional diffusion equations}. Journal of \;Differential Equations, Vol. 199(2): 211–255, 2004.

\bibitem {Fri83} A. Friedman.\quad \textit{Partial Differential Equations of Parabolic Type}.\quad Robert E. Krieger Publishing Company, Florida, 1983. Reimpresión de la original de 1964.

\bibitem {GNZ20} T. Ghoul, V. Nguyen, H. Zaag. \textit{Construction of type I blowup solutions for a higher order semilinear parabolic equation}. Advances in Nonlinear Analysis, Vol. 9, No. 1: 388-412, 2020. 

\bibitem {HKNS06} N. Hayashi, E.I. Kaikina, P.I. Naumkin, I.A Shishmarev. \textit{Asymptotics for Dissipative Nonlinear Equations}. Springer-Verlag Berlin Heidelberg, 2006. 

\bibitem {IKO02} A. Ilyin, A. Kalashnikov y O. Oleynik. \textit{Second order linear equations of parabolic type}. Journal of Mathematical Sciences, Vol. 108, No. 4: 435-542, 2002.

\bibitem {Jac01} N. Jacob. \textit{Pseudo-differential operators and Markov processes Volume I}. Imperial College Press, 2001.

\bibitem {JK19} I. Johnston, V. Kolokoltsov. \textit{Green's function estimates for time-fractional evolution equations}.\quad Fractal and fractional, arXiv:1906.12157v1, 2019.

\bibitem {Kol00} V. Kolokoltsov. \textit{Symmetric stable laws and stable-like jump-diffusions}. Proc. London Math. Soc(3), 80(3): 725-768, 2000.

\bibitem {KolV00} V. Kolokoltsov. \textit{Semiclassical Analysis for Diffusions and Stochastic Processes}. Springer, 2000. 

\bibitem {Kol09} V. Kolokoltsov. \textit{The Lévy–Khintchine type operators with variable Lipschitz continuous coefficients generate linear or nonlinear Markov processes and semigroups}. Probability Theory and Related Fields, arXiv:0911.5688v1, 2009.

\bibitem {Kol11} V. Kolokoltsov. \textit{Markov Processes, Semigroups and Generators}. Walter de Gruyter GmbH and Co. KG, Berlin/New York, 2011.

\bibitem {Kol19} V. Kolokoltsov. \textit{Differential Equations on Measures and Functional Spaces}. Birkhäuser Advanced Texts Basler Lehrbücher, e-book, 2019.

\bibitem {KSVZ16} J. Kemppainen, J. Siljander, V. Vergara, R. Zacher. \textit{Decay estimates for time-fractional and other non-local in time subdiffusion equation in $\mathbb{R}^d$}.  Math. Ann. 366: 941–979, 2016.

\bibitem {KeSZ17} J. Kemppainen, J. Siljander, R. Zacher. \textit{Representation of solutions and large-time behavior for fully nonlocal diffusion equations}. Journal Differential Equations, Vol. 263: 149–201, 2017.

\bibitem {KV14} V. Kolokoltsov and M. Veretennikova. \textit{Well-posedness and regularity of the Cauchy problem for nonlinear fractional in time and space equations}. Fractional Differential Calculus, Vol. 4, No. 1: 1-30, 2014.

\bibitem {MG00} F. Mainardi, R. Gorenflo. \textit{On Mittag-Leffler-type functions in fractional evolution processes}. Journal of Computational and Applied Mathematics. Vol. 118: 283-299, 2000.

\bibitem {MJCB14} R. Metzler, J. Jeon, A. Cherstvy, E. Barkai. \textit{Anomalous diffusion models and their \;properties: non-stationarity, non-ergodicity, and ageing at the centenary of single particle tracking}. Physical Chemistry Chemical Physics, Vol. 16: 24128-24164, 2014.

\bibitem {Pru93} J. Prüss. \textit{Evolutionary Integral Equations and Applications}. Birkhäuser Verlag, Switzerland, 1993.

\bibitem {Rud02} W. Rudin. \textit{Análisis Funcional}. Editorial Reverté S.A., 2002.

\bibitem {Schi01} R. L. Schilling. \textit{Dirichlet operators and the Positive Maximum Principle}. Integral Equations and Operator Theory 41: 74-92, Birkhäuser Verlag, Basel, 2001.

\bibitem {UZ99} V. Uchaikin, V. Zolotarev. \textit{Chance and stability: stable distributions and their \;applications}. \;Monographs Modern Probability and Statistics, 1999.

\bibitem {WZ77} R. Wheeden and A. Zygmund. \textit{Measure and Integral: An Introduction to Real Analysis}. Marcel Dekker, Inc., New York, 1977. 

\bibitem {Zac08} R. Zacher. \textit{Boundedness of weak solutions to evolutionary partial integro-differential
equations with discontinuous coefficients}. Journal of Mathematical Analysis and Applications 348: 137–149, 2008.

\bibitem {Zol86} V. Zolotarev. \textit{One-dimensional Stable Distributions}. Moscow, Nauka, 1983 (in Russian). Engl. transl. in vol. 65 of Translations of Mathematical Monographs AMS, Providence, Rhode Island, 1986.
 
\end{thebibliography}

\end{document}